\documentclass[12pt]{article}

\usepackage{amssymb,amsmath,amsfonts,amssymb}
\usepackage{graphics,graphicx,color}

\def\@abssec#1{\vspace{.05in}\footnotesize \parindent .2in
{\bf #1. }\ignorespaces}

\graphicspath{{/EPSF/}{Figures/}}
\DeclareGraphicsExtensions{.eps}

\newtheorem{theorem}{Theorem}[section]
\newtheorem{lemma}[theorem]{Lemma}
\newtheorem{proposition}[theorem]{Proposition}
\newtheorem{corollary}[theorem]{Corollary}

\def \Rm {\mathbb R}
\def \Nm {\mathbb N}
\def \Cm {\mathbb C}

\def \Sm {\mathbb S}
\newcommand{\eps}{\varepsilon}

\newcommand{\dsum}{\displaystyle\sum}
\newcommand{\dint}{\displaystyle\int}

\newcommand{\conductivity}{\gamma}    		
\newcommand{\conductivitypot}{\sigma}  		
\newcommand{\diffusion}{\gamma}         		
\newcommand{\absorption}{\sigma}            	
\newcommand{\lame}{\gamma}
\newcommand{\density}{\rho}

\newcommand{\pdr}[2]{\dfrac{\partial{#1}}{\partial{#2}}}
\newcommand{\pdrr}[2]{\dfrac{\partial^2{#1}}{\partial{#2}^2}}

\newcommand{\aver}[1]{\langle {#1} \rangle}

\newcommand{\be}{\mathbf e} 
\newcommand{\bk}{\mathbf k}
\newcommand{\bl}{\mathbf l}

\newcommand{\mH}{\mathcal H}

\newcommand{\mk}{\mathfrak k}

\newcommand{\fc}{{\mathfrak c}}

\newcommand{\mg}{{\mathfrak g}}
\newcommand{\mh}{{\mathfrak h}}

\newcommand{\fs}{{\mathfrak s}}

\newcommand{\cout}[1]{}

\newcommand{\cO}{{\mathcal O}}

\newcommand{\dX}{\partial X}
\newcommand{\X}{\tilde X}

 \renewcommand{\arraystretch}{1.5}

\title{Hybrid inverse problems and internal functionals}

\author{Guillaume Bal \thanks{Department of Applied Physics and 
        Applied Mathematics, Columbia University, 
        New York NY, 10027; gb2030@columbia.edu}}

\begin{document}
 
\maketitle

\begin{abstract}
This paper reviews recent results on hybrid inverse problems, which are also called coupled-physics inverse problems of multi-wave inverse problems. Inverse problems tend to be most useful in, e.g., medical and geophysical imaging, when they combine high contrast with high resolution.  In some settings, a single modality displays either high contrast or high resolution but not both. In favorable situations, physical effects couple one modality with high contrast with another modality with high resolution. The mathematical analysis of such couplings forms the class of hybrid inverse problems. 
  
Hybrid inverse problems typically involve two steps. In a first step, a well-posed problem involving the high-resolution low-contrast modality is solved from knowledge of boundary measurements. In a second step, a quantitative reconstruction of the parameters of interest is performed from knowledge of the point-wise, internal, functionals of the parameters reconstructed during the first step. This paper reviews mathematical techniques that have been developed in recent years to address the second step.

Mathematically, many hybrid inverse problems find interpretations in terms of linear and nonlinear (systems of) equations. In the analysis of such equations, one often needs to verify that qualitative properties of solutions to elliptic linear equations are satisfied, for instance the absence of any critical points. This paper reviews several methods to prove that such qualitative properties hold, including the method based on the construction of complex geometric optics solutions.

\end{abstract}

\tableofcontents


\renewcommand{\thefootnote}{\fnsymbol{footnote}}
\renewcommand{\thefootnote}{\arabic{footnote}}

\renewcommand{\arraystretch}{1.1}





\section{Introduction}

The success of most medical imaging modalities rests on their high, typically sub-millimeter, resolution. Computerized Tomography (CT), Magnetic Resonance Imaging (MRI), or Ultrasound Imaging (UI) are typical examples of such modalities. In some situations, these modalities fail to exhibit a sufficient contrast between different types of tissues, whereas other modalities, for example based on the optical, elastic, or electrical properties of these tissues, do display such high contrast. Unfortunately, the latter modalities, such as e.g., Optical Tomography (OT), Electrical Impedance Tomography (EIT) or Elastographic Imaging (EI), involve a highly smoothing measurement operator and are thus typically low-resolution as stand-alone modalities.

Hybrid inverse problems concern the combination of a high contrast modality with a high resolution modality. By combination, we mean the existence of a physical mechanism that couples these two modalities. Several examples of physical couplings are reviewed in section \ref{sec:physmod}. A different strategy, consisting of fusing data acquired independently for two or more imaging modalities, is referred to as multi-modality imaging and is not considered in this paper.  Examples of possible physical couplings include: optics or electromagnetism with ultrasound  in Photo-Acoustic Tomography (PAT), Thermo-Acoustic Tomography (TAT) and in Ultrasound Modulated Optical Tomography (UMOT), also called Acousto-Optic Tomography (AOT); electrical currents with ultrasound in Ultrasound Modulated Electrical Impedance Tomography (UMEIT), also called Electro-Acoustic Tomography (EAT); electrical currents with magnetic resonance in Magnetic Resonance EIT (MREIT) or Current Density Impedance Imaging (CDII); and elasticity with ultrasound in Transient Elastography (TE). Some hybrid modalities have been explored experimentally whereas other hybrid modalities have not been tested yet. Some have received quite a bit of mathematical attention whereas other ones are less well understood. While more references will be given throughout the review, we refer the reader at this point to the recent books \cite{A-Sp-08,S-Sp-2011,WW-W-07} and their references for general information about practical and theoretical aspects of medical imaging.

\medskip

Reconstructions in hybrid inverse problems typically involve two steps. In a first step, an inverse problem involving the high-resolution-low-contrast modality needs to be solved. In PAT and TAT for instance, this corresponds to reconstructing the initial condition of a wave equation from available boundary measurements. In UMEIT and UMOT, this corresponds in an idealized setting to inverting a Fourier transform that is reminiscent of the reconstructions performed in MRI. In Transient Elastography, this essentially corresponds to solving an inverse scattering problem in a time-dependent wave equation. In this review, we assume that this first step has been performed.

Our interest is in the second step of the procedure, which consists of reconstructing the coefficients that display high contrasts from the mappings obtained during the first step. These mappings involve internal functionals of the coefficients of interest. Typically, if $\conductivity$ is a coefficient of interest and $u$ is the solution to a partial differential equation involving $\conductivity$, then the internal ``measurements'' obtained in the first step take the form $H(x)=\conductivity(x) u^j(x)$ for $j=1,2$ or $H(x)=\conductivity(x) |\nabla u|^j(x)$ again for $j=1,2$.

Several questions can then be raised: are the coefficients, e.g. $\conductivity$, uniquely characterized by the internal measurements $H(x)$? How stable are the reconstructions? If specific boundary conditions are prescribed at the boundary of the domain of interest, how do the answers to the above questions depend on such boundary conditions? The answers to these questions depend on the physical model of interest. However, there are important common features that we would like to present in this review.

One such feature relates to the stability of the reconstructions. Loosely speaking, an inverse problem is well-posed, or at least not severely ill-posed, when singularities in the coefficients of interest propagate into singularities in the available data.  The map reconstructed during step 1 provides local, point-wise, information about the coefficients. Singularities of the coefficient do not need to propagate to the domain's boundary  and  we thus expect resolution of hybrid modalities to be significantly improved compared to the stand-alone high-contrast-low-resolution modalities.  This will be verified on the examples reviewed here.

Another feature is the relationship between hybrid inverse problems and nonlinear partial differential equations. Typically, both the coefficient $\conductivity$ and the solution $u$ are unknown. However, for measurements of the form $H(x)=\conductivity(x) u^j(x)$, then $\conductivity$ in the equation for $u$ can be eliminated using the expression for $H(x)$. This results in a nonlinear equation for $u(x)$. The resulting nonlinear equations often do not display any of the standard features that are amenable to proofs of uniqueness, such as admitting a variational formulation with a strictly convex functional. The main objective is to obtain uniqueness and stability results for such equations, often in the presence of redundant (overdetermined) information.

A third feature shared by many hybrid inverse problems is that their solution strategies often require that the forward solution $u$ satisfy certain qualitative properties, such as for instance the absence of any critical point (points where $\nabla u=0$). The derivation of qualitative properties such as lower bounds for the modulus of a gradient is a difficult problem. In two dimensions of space, the fact that critical points of elliptic solutions are necessarily isolated is of great help. In higher dimension, such results no longer hold in general. A framework to obtain the requested qualitative behavior of the elliptic solutions is based on the so-called complex geometric optics (CGO) solutions. Such solutions, when they can be constructed, essentially allow us to treat the unknown coefficients as perturbations of known operators, typically the Laplace operator. Using these solutions, we can construct an open set of boundary conditions for which the requested property is guaranteed. This procedure provides a  restricted class of boundary conditions for which the solutions to the hybrid inverse problems are shown to be uniquely and stably determined by the internal measurements. From a practical point of view, these mathematical results confirm the physical intuition that the coupling of high contrast and high resolution modalities indeed provides reconstructions that are robust with respect to errors in the measurements.  

\medskip

The rest of this paper is structured as follows. Section \ref{sec:physmod} is devoted to the modeling of the hybrid inverse problems and the derivation of the internal measurements for the applications considered in this paper, namely: PAT, TAT, UMEIT, UMOT, TE, CDII. The following two sections present recent results of uniqueness and stability obtained for such hybrid inverse problems: Section \ref{sec:diffabs} focuses on internal functionals of the solution $u$ of the forward problem, whereas section \ref{sec:diff} is concerned with internal functionals of the gradient of the solution $\nabla u$. As we mentioned above, these uniqueness and stability results hinge on the forward solutions $u$ to verify some qualitative properties. Section \ref{sec:propsol} summarizes some of these properties in the two-dimensional case and presents the derivation of such properties in higher spatial dimensions by means of complex geometric optics (CGO) solutions. Some concluding remarks are proposed in section \ref{sec:conclu}.



\section{Physical modeling}
\label{sec:physmod}

 High resolution imaging modalities include Ultrasound Imaging and Magnetic Resonance Imaging. High contrast modalities include Optical Tomography, Electrical Impe\-dance Tomography, and Elastography. This sections briefly presents four couplings between high-contrast and high-resolution modalities: two different methods to couple ultrasound and optics or (low frequeny) electromagnetism in PAT/TAT via the photo-acoustic effect and in UMOT/UMEIT via ultrasound modulation; the coupling between Ultrasound and Elastography in Transient Elastography; and the coupling between Electrical Impedance Tomography and Magnetic Resonance Imaging in CDII/MREIT.

\subsection{The Photo-acoustic effect}
\label{sec:photoacoustic}

The photoacoustic effect may be described as follows. A pulse of radiation is sent into a domain of interest. A fraction of the propagating radiation is absorbed by the medium. This generates a thermal expansion, which is the source of ultrasonic waves. Ultrasound then propagates to the boundary of the domain where ultrasonic transducers measure the pressure field. The physical coupling between the absorbed radiation and the emitted sound is called the photoacoustic effect. This is the premise for the medical imaging technique Photoacoustic Tomography (PAT). 

Two types of radiation are typically considered. In Optoacoustic Tomography (OAT), near-infra-red photons, with wavelengths typically between $600nm$ and $900nm$ are used. The reason for this frequency window is that they are not significantly absorbed by water molecules and thus can propagate relatively deep into tissues. OAT is often simply referred to as PAT and we will follow this convention here. In Thermoacoustic Tomography (TAT), low frequency microwaves, with wavelengths on the order of $1m$, are sent into the medium. The rationale for using such frequencies is that they are less absorbed than optical frequencies and thus propagate into deeper tissues.

In both PAT and TAT, the first step of an inversion procedure is the reconstruction of the map of absorbed radiation from the ultrasonic measurements. In both applications, the inversion may be recast as the reconstruction of an initial condition of a wave equation from knowledge of ultrasound measurements. Assuming a domain of infinite extension with non-perturbative measurements to simplify the presentation, ultrasound propagation is modeled by the following wave equation
\begin{equation}
  \label{eq:wave}
  \begin{array}{ll}
  \dfrac{1}{c_s^2(x)} \pdrr pt - \Delta p =0,\quad &t>0,\,\,x\in \Rm^n \\
  p(0,x) = H(x) \quad\mbox{ and } \quad 
  \pdr pt(0,x) =0 \quad &x\in\Rm^n.
  \end{array}
\end{equation}
Here $c_s$ is the sound speed assumed to be known, $n$ is spatial dimension, and $H(x)$ is the ultrasonic signal generated at time $t=0$. Measurements are then of the form $p(t,x)$ for $t>0$ and $x\in\partial X$ at the boundary of a domain $X$ where $H(x)$ is supported.

Note that the effect of propagating radiation is modeled as an {\em initial} condition at $t=0$. The reason for this stems from the large difference between light speed (roughly $2.3\,10^8m/s$ in water) and sound speed (roughly $1.5\,10^3m/s$ in water). When a short pulse of radiation is emitted into the medium, we may assume that it propagates into the medium at a time scale that is very short compared to that of ultrasound.  This is a very valid approximation in PAT but is a limiting factor in the (still significantly sub-millimeter) spatial resolution we expect to obtain in TAT; see, e.g., \cite{BJJ-IP-10,BRUZ-IP-11}.

For additional references to the photo-acoustic effect, we refer the reader to the works \cite{CAB-JOSA-09,CLB-SPIE-09,FSS-PRE-07,XW-RSI-06,XWKA-CRC-09} and their references. The first step in thermo- and photo-acoustics is the reconstruction of the absorbed radiation map $H(x)$ from boundary acoustic wave measurements. There is a vast literature on this inverse source problem in the mathematical and physical literatures.  We refer the reader to e.g. \cite{ABJK-SR-09,FPR-JMA-04,HSBP-IP-04,HKN-IP-08,KK-EJAM-08,PS-IP-07,SU-IP-09}. Serious difficulties may need to be addressed in this first step, such as e.g. limited data,  spatially varying acoustic sound speed \cite{ABJK-SR-09,HKN-IP-08,SU-IP-09}, and the  effects of acoustic wave attenuation \cite{KowSch10}. In this paper, we assume that the absorbed radiation map $H(x)$ has been reconstructed. This provides now {\em internal} information about the properties of the domain of interest. What we can extract from such information depends on the model of radiation propagation. The resulting inverse problems are called Quantitative PAT (QPAT) and Quantitative TAT (QTAT) for the different modalities of radiation propagation, respectively.

In the PAT setting with near-infra-red photons, arguably the most accurate model for radiation propagation is the radiative transfer equation. We shall not describe this model here and refer the reader to \cite{BJJ-IP-10} for QPAT in this setting and to \cite{B-IP-09} for more general inverse problems for the radiative transfer equation. The models we consider for radiation propagation are as follows.

\subsubsection{QPAT modeling}
\label{sec:QPATmodel}
In the diffusive regime, photon (radiation) propagation is modeled by the following second-order elliptic equation
\begin{equation}\label{eq:diff}
	\begin{array}{ll}
		-\nabla\cdot \diffusion(x)\nabla u + \absorption(x) u = 0 \quad  & \text{in }\ X\\
		u=f &\text{on }\ \dX.
	\end{array}
\end{equation}
To simplify, we assume that Dirichlet conditions are prescribed at the boundary of the domain $\dX$. Throughout the paper, we assume that $X$ is a bounded open domain in $\Rm^n$ with smooth boundary $\dX$. The optical coefficients $(\diffusion(x),\absorption(x))$  are  $\diffusion(x)$ the diffusion coefficient and $\absorption(x)$ the absorption coefficient, which are assumed to be bounded from above and below by positive constants.

The information about the coefficients in QPAT takes the following form:
\begin{equation}
  \label{eq:HQPAT}
    H(x) = \Gamma(x) \absorption(x) u(x) \quad \mbox{ a.e. } x\in X.
\end{equation}
The coefficient $\Gamma(x)$ is the Gr\"uneisen coefficient. It models the strength of the photo-acoustic effect, which converts absorption of radiation into emission of ultrasound. The objective of QPAT is to reconstruct $(\diffusion,\absorption,\Gamma)$ from knowledge of $H(x)$ in \eqref{eq:HQPAT} obtained for a given number of illuminations $f$ in \eqref{eq:diff}.  This is an example of an internal measurement that is {\em linear} in the solution $u(x)$ and the absorption coefficient $\absorption$. For references on QPAT, see, e.g., \cite{BR-OL-11,BU-IP-10,CAB-JOSA-09,CLB-SPIE-09,RN-PRE-05,Z-AO-10} and their references. 

\subsubsection{QTAT modeling}
\label{sec:QTATmodel}

Low frequency Radiation in QTAT is modeled by the following system of Maxwell's equations:
\begin{equation}
  \label{eq:invbrdyE}
  \begin{array}{ll}
   -\nabla\!\times\nabla\!\times E + k^2 E + i k \conductivitypot(x) E =0, \qquad& X \\
   \nu\times E=f \qquad &\partial X.
  \end{array}
\end{equation}
Here, $E$ is the (time-harmonic) electromagnetic field with fixed wavenumber $k=\frac\omega c$ where $\omega$ is frequency and $c$ light speed.  We assume that radiation is controlled by the boundary condition $f(x)$ on $\dX$.  The unknown coefficient is the conductivity (absorption) coefficient $\conductivitypot(x)$.
Setting $\Gamma=1$ to simplify, the map of  absorbed electromagnetic radiation is then of the form
\begin{equation}
  \label{eq:IScIDE}
  H(x) = \conductivitypot(x) |E|^2(x).
\end{equation}

The above system of equations may be simplified by modeling radiation by a scalar quantity $u(x)$. In this setting, radiation is modeled by the following Helmholtz equation:
\begin{equation}
  \label{eq:invbrdy}
  \begin{array}{ll}
   \Delta u + k^2 u + i k \conductivitypot(x) u =0, \qquad& X \\
   u=f \qquad &\partial X,
  \end{array}
\end{equation}
for a given boundary condition $f(x)$.  The internal data are then of
the form
\begin{equation}
  \label{eq:IScID}
  H(x) = \conductivitypot(x) |u|^2(x).
\end{equation}
For such models, QTAT then consists of reconstructing $\conductivitypot(x)$ from knowledge of $H(x)$. Note that $H(x)$ is now a {\em quadratic} quantity in the solutions $E(x)$ or $u(x)$. There are relatively few results on QTAT; see \cite{BRUZ-IP-11,LPKW-PRE-08}.

\subsection{The ultrasound modulation effect}
\label{sec:UM}

We consider the following elliptic equation
\begin{equation}
  \label{eq:ellipticDsigma}
  -\nabla\cdot \diffusion(x)\nabla u + \absorption(x) u = 0\quad \mbox{ in } X, \qquad u=f\quad \mbox{ on } \partial X.
\end{equation}
The objective of ultrasound modulation is to send an acoustic signal through the domain $X$ that modifies the coefficients $\diffusion$ and $\absorption$. We assume here that the sound speed is constant and that we are able to generate an acoustic signal that takes the form of the plane wave $p\cos(k\cdot x + \varphi)$ where $p$ is amplitude, $k$ wave-number and $\varphi$ an additional phase. We assume that the acoustic signal modifies the properties of the diffusion equation and that the effect is small. The coefficients in \eqref{eq:ellipticDsigma} are thus modified as 
\begin{equation}
  \label{eq:modifcoefs}
  \diffusion_\eps(x) = \diffusion(x) (1+\zeta \eps \fc ) + O(\eps^2),\qquad
  \absorption_\eps(x) = \absorption(x) (1+\eta \eps \fc ) + O(\eps^2),
\end{equation}
where we have defined $\fc=\fc(x)=\cos(k\cdot x + \varphi)$ and where $\eps=p\Gamma$ is the product of the acoustic amplitude $p\in\Rm$ and a measure $\Gamma>0$ of the coupling between the acoustic signal and the modulations of the constitutive parameters in \eqref{eq:ellipticDsigma}. We assume that $\eps\ll1$. The terms in the expansion are characterized by $\zeta$ and $\eta$ and depend on the specific application. 

Let $u$ and $v$ be solutions of \eqref{eq:ellipticDsigma} with fixed boundary conditions $f$ and $h$, respectively. When the acoustic field is turned on, the coefficients are modified as described in \eqref{eq:modifcoefs} and we denote by $u_\eps$ and $v_\eps$ the corresponding solution. Note that $u_{-\eps}$ is the solution obtained by changing the sign of $p$ or equivalently by replacing $\varphi$ by $\varphi+\pi$.

By standard regular perturbation arguments, we find that $u_\eps=u_0+\eps u_{1} + O(\eps^2)$.
Multiplying the equation for $u_\eps$ by $v_{-\eps}$ and the equation for $v_{-\eps}$ by $u_\eps$, subtracting the results, and using standard integrations by parts, we obtain that 
\begin{equation}
  \label{eq:Green}
  \dint_X (\diffusion_\eps-\diffusion_{-\eps})\nabla u_\eps\cdot\nabla v_{-\eps} + (\absorption_\eps-\absorption_{-\eps}) u_\eps v_{-\eps} dx = \dint_{\partial X}  \diffusion_{-\eps} \pdr{v_{-\eps}}{\nu} u_\eps -\diffusion_\eps \pdr{u_\eps}{\nu} v_{-\eps}  d\sigma.
\end{equation}
We assume that $\diffusion_\eps \partial_\nu u_\eps$ and $\diffusion_\eps \partial_\nu v_\eps$ are measured on $\partial X$, at least on the support of $v_{\eps}=h$ and $u_\eps=f$, respectively, for several values $\eps$ of interest. The above equation also holds if the Dirichlet boundary conditions are replaced by Neumann boundary conditions. Let us define
\begin{equation}
  \label{eq:measbdry}
  J_\eps :=  \dfrac12\dint_{\partial X}  \diffusion_{-\eps} \pdr{v_{-\eps}}{\nu} u_\eps -\diffusion_\eps \pdr{u_\eps}{\nu} v_{-\eps}  d\sigma  \,\,=\,\, \eps J_1 + \eps^2 J_2 + O(\eps^3).
\end{equation}
We assume that the real valued functions $J_m=J_m(k,\varphi)$ are known from the physical measurement of the Cauchy data of the form $(u_\eps,\diffusion_\eps \partial_\nu u_\eps)$ and  $(v_\eps,\diffusion_\eps \partial_\nu v_\eps)$ on $\partial X$. 

Equating like powers of $\eps$, we find that at the leading order
\begin{equation}
  \label{eq:leadingorder}
  \dint_X \big[\zeta \diffusion(x) \nabla u_0\cdot\nabla v_0(x)+ \eta \absorption(x) u_0v_0(x) \big] \cos(k\cdot x+\varphi) dx = J_1(k,\varphi).
\end{equation}
Acquiring this for all $k\in\Rm^n$ and $\varphi=0,\frac\pi2$, this yields after inverse Fourier transform:
\begin{equation}
  \label{eq:H0Lap}
  H[u_0,v_0](x) = \zeta \diffusion(x) \nabla u_0\cdot\nabla v_0(x) + \eta \absorption(x) u_0v_0(x).
\end{equation}


In the setting of ultrasound modulated optical tomography (UMOT), the coefficients $\diffusion_\eps$ and $\absorption_\eps$ in \eqref{eq:modifcoefs} take the form \cite{BS-PRL-10}
\begin{math} 
   \diffusion_\eps(x) = \frac{\tilde \diffusion_\eps}{c_\eps^{n-1}}(x)\mbox{ and }
   \absorption_\eps(x) = \frac{\tilde \absorption_\eps}{c_\eps^{n-1}}(x),
 \end{math}
where $\tilde\absorption_\eps$ is the absorption coefficient, $\tilde \diffusion_\eps$ is the diffusion coefficient, $c_\eps$ is the light speed, and $n$ is spatial dimension. When the pressure field is turned on, the amount of scatterers and absorbers is modified by compression and dilation. Since the diffusion coefficient is inversely proportional to the scattering coefficient, we find that 
\begin{displaymath} 
   \tilde\absorption_\eps(x) = \tilde\absorption(x) \big( 1+\eps \fc(x)\big),\qquad \dfrac1{\diffusion_\eps(x)}=\dfrac1{\diffusion(x)}\big(1+\eps \fc(x)\big).
 \end{displaymath}
The pressure field changes the index of refraction of light as follows
\begin{math} 
    c_\eps(x) = c(x)(1+\psi\eps \fc(x)),
\end{math}
where $\psi$ is a constant (roughly equal to $\frac13$ for water). This shows that
\begin{equation}
  \label{eq:zetaeta}
  \zeta = - (1+(n-1)\psi),\qquad \eta = 1-(n-1)\psi.
\end{equation}

In the application of ultrasound modulated electrical impedance tomography (UMEIT), $\diffusion(x)$ is a conductivity coefficient and $\absorption=0$. We then have $\diffusion_\eps(x)=\diffusion(x)(1+\eps \fc(x))$ with thus $\zeta=1$ and $\eta=0$.  The objective of UMOT and UMEIT is to reconstruct (part of) the coefficients $(\diffusion(x),\absorption(x))$ in the elliptic equation
\begin{equation}
  \label{eq:ellipticDsigma1}
  -\nabla\cdot \diffusion(x)\nabla u + \absorption(x) u = 0\quad \mbox{ in } X, \qquad u=f\quad \mbox{ on } \partial X.
\end{equation}
from measurements of the form 
\begin{equation}
  \label{eq:H0Lap1}
  H[u_0,v_0](x) = \zeta \diffusion(x) \nabla u_0\cdot\nabla v_0(x) + \eta \absorption(x) u_0v_0(x),
\end{equation}
for one or several values of the illumination $f(x)$ on $\dX$.

In a simplified version of UMOT (also called acousto-optic tomography; AOT), $\zeta=0$ and the measurements are quadratic (or bilinear) in the solutions to the elliptic equation. More challenging mathematically is the case $\zeta=1$ and $\eta=0$ where the measurements are quadratic (or bilinear) in the {\em gradients} of the solution. No theoretical results exist to date in the setting where both $\zeta$ and $\eta$ are non-vanishing.

The effect of ultrasound modulation is difficult to observe experimentally as the coupling coefficient $\Gamma$ above is rather small. For references on ultrasound modulation in different contexts, we refer the reader to, e.g., \cite{ABCTF-SIAP-08,B-UMEIT-11,BS-PRL-10,CFGK-SJIS-09,GS-SIAP-09,KK-AET-11,ZW-SPIE-04}. These references concern the so-called incoherent regime of wave propagation, while the coherent regime, whose mathematical structure is different, is addressed in the physical literature in, e.g., \cite{AFRBG-OL-05,KLZZ-JOSA-97,W-JDM-04}.

%
%


\subsection{Transient Elastography}
\label{sec:elast}

Transient elastography images the (slow) propagation of shear waves using ultrasound. For more details, see, e.g., \cite{MZM-IP-10} and its extended list of references. As shear waves propagate, the resulting displacements can be imaged by ultra-fast ultrasound. Consider a scalar approximation of the equations of elasticity
\begin{equation}
  \label{eq:elastscalar0}
  \begin{array}{ll}
  \nabla\cdot \lame(x)\nabla u(x, t) = \density(x)\partial_{tt} u (x, t),\quad & t\in\Rm,\,x\in X\\
  u(x,t) = f(x,t),\quad & t\in\Rm,\, x\in\dX,
  \end{array}
\end{equation}
where $u(x,t)$ is the (say, downward) displacement, $\lame(x)$ is one of the Lam\'e parameters and $\density(x)$ is density. Using ultra-fast ultrasound measurements, the displacement $u(x,t)$ can be imaged. This results in a very simplified model of transient elastography where we aim to reconstruct $(\lame,\density)$ from knowledge of $u(x,t)$; see \cite{MZM-IP-10} for more complex models. We may slightly generalize the model as follows. Upon taking Fourier transforms in the time domain and accounting for possible dispersive effects of the tissues, we obtain
\begin{equation}
  \label{eq:elastscalar}
  \begin{array}{ll}
  \nabla\cdot \lame(x;\omega)\nabla u(x;\omega) +\omega^2 \density(x;\omega) u (x;\omega)=0,\quad & \omega\in\Rm,\,x\in X\\
  u(x;\omega) = f(x;\omega),\quad & \omega\in\Rm,\, x\in\dX.
  \end{array}
\end{equation}
The inverse transient elastography problem with dispersion effect would then be the reconstruction of $(\lame(x;\omega),\density(x;\omega))$ from knowledge of $u(x;\omega)$ corresponding to one or several boundary conditions $f(x;\omega)$ applied at the boundary $\dX$. This hybrid inverse problem again involves measurements that are {\em linear} in the solution $u$.

\subsection{Current Density Imaging}

Magnetic Impedance Electrical Impedance Tomography (MREIT) and Current Density Impedance Imaging (CDII) are two modalities aiming to reconstruct the conductivity in an equation using magnetic resonance imaging (MRI). The electrical potential $u$ solves the following elliptic equation
\begin{equation}
  \label{eq:ellipticsigma}
  -\nabla\cdot \conductivity(x)\nabla u = 0\quad \mbox{ in } X, \qquad u=f\quad \mbox{ on } \partial X,
\end{equation}
with $\conductivity(x)$ the unknown conductivity and $f$ a prescribed voltage at the domain's boundary. The electrical current density $J=-\conductivity\nabla u$ satisfies the system of Maxwell's equations
\begin{equation}
  \label{eq:MaxwellJ}
  \nabla\cdot J =0,\qquad J = \dfrac{1}{\mu_0} \nabla\times B, \qquad x\in X.
\end{equation}
Here $\mu_0$ is a constant, known, magnetic permeability.

Ideally, the whole field $B$ can be reconstructed from MRI measurements. This provides access to the current density $J(x)$ in the whole domain $X$. CDI then corresponds to reconstructing $\conductivity$ from knowledge of $J$. In practice, acquiring $B$ requires rotation of the domain of interest (or of the MRI apparatus) which is not straightforward. MREIT thus assumes knowledge of the third component $B_z$ of the magnetic field for several possible boundary conditions. This provides information about $\conductivity(x)$. We do not consider the MREIT inverse problem further and refer the reader to the recent review \cite{SW-SR-11} and its references for additional information. 

Several works have considered the problem of the reconstruction of $\conductivity$ in \eqref{eq:ellipticsigma}  from knowledge of the scalar information $|J|$ rather than the full current $J$. This inverse problem, referred to as the $1-$Laplacian, will be addressed below and compared to the $0-$Laplacian that appears in UMEIT and UMOT. For references on MREIT and CDII, we refer the reader to \cite{KKSY-SIMA-02,NTT-IP-07,NTT-IP-09,NTT-Rev-11} and their references.

\section{Reconstructions from functionals of $u$}
\label{sec:diffabs}

In this section, we consider internal measurements $H(x)$ of the form $H(x)=\tau(x)u(x)$ for $\tau(x)$ a function that depends linearly on unknown coefficients such as the diffusion coefficient $\diffusion$ or the absorption coefficient $\absorption$ in section \ref{sec:PAT} and internal measurements $H(x)$ of the form $H(x)=\tau(x)|u(x)|^2$ in section \ref{sec:TAT}, where $\tau$ again depends linearly on unknown coefficients. Measurements of the first form find applications in Quantitative Photo-acoustic Tomography (QPAT) and Transient Elastography (TE) while measurements of the second form find applications in Quantitative Thermo-Acoustic Tomography (QTAT) and simplified models of Acousto-Optics Tomography (AOT).


\subsection{Reconstructions from linear functionals in $u$}
\label{sec:PAT}

Recall the elliptic model for photon propagation in tissues:
\begin{equation}
  \label{eq:ellipticDsigma2}
  -\nabla\cdot \diffusion(x)\nabla u + \absorption(x) u = 0\quad \mbox{ in } X, \qquad u=f\quad \mbox{ on } \partial X.
\end{equation}
The information about the coefficients in QPAT takes the following form:
\begin{equation}
  \label{eq:HQPAT2}
    H(x) = \Gamma(x) \absorption(x)  u(x) \quad \mbox{ a.e. } x\in X.
\end{equation}
The coefficient $\Gamma(x)$ is the Gr\"uneisen coefficient. In many works in QPAT, it is assumed to be constant. We assume here that it is Lipschitz continuous and bounded above and below by positive constants.

\subsubsection{Non-unique reconstruction of three coefficients.}
\label{sec:recthreecoefs}

Let $f_1$ and $f_2$ be two Dirichlet conditions on $\partial X$ and $u_1$ and $u_2$ be the corresponding solutions to \eqref{eq:ellipticDsigma2}. We make the following assumptions:

\begin{itemize}
\vspace{-.3cm}
\item[(i)] The coefficients $(\diffusion,\absorption,\Gamma)$ are of class $W^{1,\infty}(X)$ and bounded above and below by positive constants. The coefficients $(\diffusion,\absorption,\Gamma)$ are known on $\dX$.
\vspace{-.3cm}
\item[(ii)]
The  illuminations $f_1$ and $f_2$ are positive functions on $\dX$ and are the traces on $\dX$ of functions of class $C^3(\bar X)$.
\vspace{-.3cm}
\item[(iii)] the vector field 
\begin{equation}
  \label{eq:beta}
  \beta := H_1\nabla H_2 -  H_2\nabla H_1 = H_1^2 \nabla \frac{H_2}{H_1} = H_1^2\nabla \frac{u_2}{u_1} = - H_2^2 \nabla \frac{H_1}{H_2}
\end{equation} 
is a vector field in $W^{1,\infty}(X)$ such that $\beta\not\equiv0$ (on a set of positive measure).
\vspace{-.3cm}
\item[(iii')] same as (iii) above with 
\begin{equation}
  \label{eq:lowerbd}
  |\beta|(x) \geq \alpha_0>0 ,\qquad \mbox{ a.e. } x\in \X.
\end{equation} 
\end{itemize}
Beyond the regularity assumptions on $(\diffusion,\absorption,\Gamma)$, the domain $X$, and the boundary conditions $f_1$ and $f_2$, the only real assumption we impose is \eqref{eq:lowerbd}. In general, there is no guaranty that the gradient of $\frac{u_2}{u_1}$ does not vanish. 
In dimension $d=2$, a simple condition guarantees that \eqref{eq:lowerbd} holds. We have the following result \cite{A-AMPA-86,NTT-IP-07}:
\begin{lemma}[\cite{BR-IP-11}]\label{lem:crit}
    Assume that $h=\frac{g_2}{g_1}$ on $\dX$ is an almost two-to-one function in the sense of \cite{NTT-IP-07}, i.e., a function that is a two-to-one map except possibly at its minimum and at its maximum. Then \eqref{eq:lowerbd} is satisfied. 
\end{lemma}

In dimension $d\geq3$, the above result on the (absence of) critical points of elliptic solutions no longer holds. By continuity, we verify that  \eqref{eq:lowerbd} is satisfied for a large class of illuminations when $\diffusion$ is close to a constant and $\absorption$ is sufficiently small. 
For arbitrary coefficients $(\diffusion,\absorption)$ in dimension $d\geq3$, a proof based on CGO solutions shows that \eqref{eq:lowerbd} is satisfied  for an open set of illuminations; see \cite{BU-IP-10} and section \ref{sec:CGO} below. Note also that \eqref{eq:lowerbd} is a sufficient condition for us to solve the inverse problem of QPAT. In \cite{A-AMPA-86}, a similar problem is addressed in dimension $d=2$ without assuming a constraint of the form \eqref{eq:lowerbd}. 

We first prove a  result that provides uniqueness up to a specified transformation.
\begin{theorem}[\cite{BR-IP-11,BU-IP-10}]
  \label{thm:unique} Assume that hypotheses {\em (i)-(iii)} hold.  Then 
  \begin{itemize}
  \vspace{-.25cm}\item[(a)]$H_1(x)$ and $H_2(x)$ uniquely determine the measurement operator ${\cal H}: H^{\frac12}(\dX)\to H^1(X)$, which to $f$ defined on $\dX$ associates ${\cal H}(f) = H$ in $X$ defined by \eqref{eq:HQPAT}.
  \vspace{-.25cm}\item[(b)] The measurement operator ${\cal H}$ uniquely determines the two functionals:
  \begin{equation}
  \label{eq:chiq}
    \chi(x) := \dfrac{\sqrt \diffusion}{\Gamma\absorption}(x), \qquad q(x) := -\Big(\dfrac{\Delta \sqrt \diffusion}{\sqrt \diffusion} + \dfrac{\absorption}\diffusion \Big)(x).
\end{equation}
Here $\Delta$ is the Laplace operator.
\vspace{-.25cm}\item[(c)] Knowledge of the two functionals $\chi$ and $q$ uniquely determines $H_1(x)$ and $H_2(x)$. In other words, the reconstruction of $(\diffusion,\absorption,\Gamma)$ is unique up to transformations that leave $(\chi,q)$ invariant.
\end{itemize}
\end{theorem}
The proof of this theorem is given in \cite{BR-IP-11} under the additional assumption (iii'). The following minor modification allows one to prove the theorem as stated above. The proof in \cite{BR-IP-11} is based on the fact that $\int_X(\rho-1)^2|\beta|^2 dx=0$ implies that $\rho=1$ a.e. Under assumption (iii), $\beta$ does not vanish on a set of positive (Lebesgue) measure. However, $\beta$ is the solution of \eqref{eq:trchi} below with $\chi^2$ bounded from below by a positive constant. This implies that $\frac{u_2}{u_1}$ is the solution of an elliptic equation. If $\beta$ vanishes on a ball, and hence $u_2=Cu_1$ for some constant $C$, on a set of positive measure, hence on a ball, then $u_2-Cu_1$ is a vanishing solution of \eqref{eq:ellipticDsigma2}, and hence vanishes everywhere by the unique continuation principle; see e.g. \cite[Chapter 3]{isakov-98}. This would violate (iii). 

\subsubsection{Reconstruction of two coefficients.}
\label{sec:rectwocoefs}

The above result shows that the unique reconstruction of $(\diffusion,\absorption,\Gamma)$ is not possible even from knowledge of the full measurement operator ${\cal H}$ defined in Theorem \ref{thm:unique}. Two well-chosen illuminations uniquely determine the functionals $(\chi,q)$ and acquiring additional measurements does not provide any new information. However, we can prove that if one coefficient in $(\diffusion,\absorption,\Gamma)$ is known, then the other two coefficients are uniquely determined:
\begin{corollary}[\cite{BR-IP-11}]\label{cor:uniq}
   Under the hypotheses of the previous theorem, let $(\chi,q)$ in \eqref{eq:chiq} be known. Then:
   \begin{itemize}
\vspace{-.25cm}\item[(a)] If $\Gamma$ is known, then $(\diffusion,\absorption)$ are uniquely determined.
\vspace{-.25cm}\item[(b)] If $\diffusion$ is known, then $(\absorption,\Gamma)$ are uniquely determined.
\vspace{-.25cm}\item[(c)] If $\absorption$ is known, then $(\diffusion,\Gamma)$ are uniquely determined.
\end{itemize}
\end{corollary}

The above uniqueness results are {\em constructive}. In all cases, we need to solve the following transport equation for $\chi$:
\begin{equation}
  \label{eq:trchi}
  -\nabla \cdot (\chi^2\beta) =0 \quad \mbox{ in } X,\qquad \chi_{|\dX} \mbox{ known on }\dX,
\end{equation}
 with $\beta$ the vector field defined in \eqref{eq:beta}. This uniquely defines $\chi>0$. Then we find that 
 \begin{equation}
  \label{eq:q}
  q(x) = -\dfrac{\Delta (H_1\chi)}{H_1\chi} = -\dfrac{\Delta (H_2\chi)}{H_2\chi}. 
\end{equation}
This provides explicit reconstructions for $(\chi,q)$ from knowledge of $(H_1,H_2)$ when \eqref{eq:lowerbd} holds.

In case (b), no further equation needs to be solved. In cases (a) and (c), we need to solve an elliptic equation for $\sqrt \diffusion$, which is the linear equation 
 \begin{equation}
  \label{eq:sqrtD}
     (\Delta+q) \sqrt{\diffusion} + \dfrac{1}{\Gamma\chi}  =0, \quad \X, \qquad \sqrt \diffusion_{|\dX} = \sqrt{\diffusion_{|\dX}}, \quad \dX,
\end{equation} 
in (a) and the (uniquely solvable) nonlinear (semi-linear) equation 
\begin{equation}\label{eq:sqrtD Semi}
    \sqrt \diffusion (\Delta + q) \sqrt \diffusion  +\absorption =0 \quad \X, \qquad \sqrt \diffusion_{|\dX} = \sqrt{\diffusion_{|\dX}}, \quad \dX,
 \end{equation} 
 in (c). These inversion formulas were implemented numerically in \cite{BR-IP-11}. Moreover, reconstructions are known to be H\"older or Lipschitz stable depending on the metric used in the stability estimate. For instance, we have:
\begin{theorem} [\cite{BR-IP-11}]
  \label{thm:stab}
  Assume that the hypotheses of Theorem \ref{thm:unique} and {\em (iii')} hold. Let $H=(H_1,H_2)$ be the measurements corresponding to the coefficients $(\diffusion,\absorption,\Gamma)$ for which hypothesis (iii) holds. Let $\tilde H=(\tilde H_1,\tilde H_2)$ be the measurements corresponding to the same illuminations $(f_1,f_2)$ with another set of coefficients $(\tilde \diffusion,\tilde \absorption,\tilde \Gamma)$ such that (i) and (ii) still hold. Then we find that 
  \begin{equation}
  \label{eq:stab1}
   \|\chi - \tilde \chi \|_{L^p(X)} \leq C \|H-\tilde H\|_{(W^{1,\frac p2}(X))^2}^{\frac12} ,\qquad \mbox{ for all } 2\leq p<\infty.
\end{equation}
Let us assume, moreover, that $\gamma(x)$ is of class $C^3(\bar X)$. Then we have that
 \begin{equation}
  \label{eq:stab3}
   \|\chi - \tilde \chi \|_{L^\infty(X)} \leq C \|H-\tilde H\|_{(L^{\frac p2}(X))^2}^{\frac p{3(d+p)}} ,\qquad \mbox{ for all } 2\leq p<\infty.
\end{equation}
\end{theorem}
We may for instance choose $p=4$ above to measure the noise level in the measurement $H$ in the square integrable norm when noise is described by its power spectrum in the Fourier domain. This shows that reconstructions in QPAT are H\"older stable, unlike the corresponding reconstructions in Optical Tomography \cite{B-IP-09,U-IP-09}.

\subsubsection{An application to Transient Elastography}

We can apply the above results to the time-harmonic reconstruction in a simplified model of transient elastography. 
Let us assume that $\lame$ and $\density$ are unknown functions of $x\in X$ and $\omega\in\Rm$. Recall that the displacement solves \eqref{eq:elastscalar}. Assuming that $u(x;\omega)$ is known after step 1 of the reconstruction using the ultrasound measurements, then we are in the setting of Theorem \ref{thm:unique} with $\Gamma\absorption=1$. Let us then assume that the two illuminations $f_1(x;\omega)$ and $f_2(x;\omega)$ are chosen such that for $u_1$ and $u_2$ the corresponding solutions of \eqref{eq:elastscalar}, we have that \eqref{eq:lowerbd} holds. We have seen a sufficient condition for this to hold in dimension $n=2$ in Lemma \ref{lem:crit} and will present other sufficient conditions in section \ref{sec:CGO} below devoted to CGO solutions in the setting $n\geq3$. Then,  \eqref{eq:chiq} shows that the reconstructed function $\chi$ uniquely determines the Lam\'e parameter $\lame(x;\omega)$ and that the reconstructed function $q$ then uniquely determines $\omega^2\density$ and hence the density parameter $\density(x;\omega)$. The reconstructions are performed for each frequency $\omega$ independently. We may summarize this as follows:
\begin{corollary}\label{cor:uniq2}
   Under the hypotheses Theorem \ref{thm:unique} and the hypotheses described above, let $(\chi,q)$ in \eqref{eq:chiq} be known. Then $(\lame(x;\omega),\density(x;\omega))$ are uniquely determined by two well-chosen measurements. Moreover, the stability results in Theorem \ref{thm:stab} hold. 
\end{corollary}

Alternatively, we may assume that in a given range of frequencies, $\lame(x)$ and $\density(x)$ are independent of $\omega$. In such a setting, we expect that one measurement $u(x;\omega)$ for two different frequencies will provide sufficient information to reconstruct $(\lame(x),\density(x))$. Assume that $u(x;\omega)$ is known for $\omega=\omega_j$, $j=1,2$ and define $0<\alpha=\omega_2^2\omega_1^{-2}\not=1$. Then straightforward calculations show that 
\begin{equation}\label{eq:betaalpha} 
     \nabla\cdot \lame \beta_\alpha =0, \quad \beta_\alpha= \big(u_1\nabla u_2-\alpha u_2\nabla u_1).
 \end{equation}
This provides a transport equation for $\lame$ that can be solved stably provided that $|\beta_\alpha|\geq c_0>0$, i.e., $\beta_\alpha$ does not vanish on $X$. Then, Theorem \ref{thm:unique} and  Theorem \ref{thm:stab} apply in this setting. Since $\beta_\alpha$ cannot be written as the ratio of two solutions as in \eqref{eq:beta} when $\alpha=1$, the results obtained in Lemma \ref{lem:crit} do not apply when $\alpha\not=1$. However, we prove in section \ref{sec:CGO} that  $|\beta_\alpha|\geq c_0>0$ is satisfied for an open set of illuminations constructed by means of CGO solutions for all $\alpha>0$; see \eqref{eq:lowerbetaalpha} below.


\subsubsection{Reconstruction of one coefficient}

Let us conclude this section by some comments on the reconstruction of a single coefficient from a measurement linear in $u$. From an algorithmic point of view, such reconstructions are siginficantly simpler. Let us consider  the framework of Corollary \ref{cor:uniq}. When $\Gamma$ is the only unknown coefficient, then we solve for $u$ in \eqref{eq:ellipticDsigma2} and reconstruct $\Gamma$ from knowledge of $H$. 

When only $\absorption$ is unknown, then we solve the elliptic equation for $u$ 
\begin{displaymath} 
   -\nabla\cdot \diffusion\nabla u + \dfrac{H}{\Gamma} =0 \quad \mbox{ in } X,\qquad u=g\quad \mbox{ on } \partial X,
\end{displaymath}
and then evaluate $\sigma=\frac{H}{\Gamma u}$.

When only $\diffusion$ is unknown with either $H=\absorption u$ in QPAT or with $H=u$ in elastography or in applications to ground water flows \cite{A-AMPA-86,R-SIAP-81}, then $u$ is known and $\diffusion$ solves the following transport equation 
\begin{displaymath} 
-\nabla\cdot \diffusion \nabla u = S \quad \mbox{ in } X,\qquad \diffusion=\diffusion_{|\partial X} \quad \mbox{ on } \partial X,
\end{displaymath}
with $S$ known. Provided that the vector field $\nabla u$ does not vanish, the above equation admits a unique solution as in \eqref{eq:trchi}. The stability results of Theorem \ref{thm:stab} then apply. Other stability results based on solving the transport equation by the method of characteristics are presented in \cite{R-SIAP-81}. In two dimensions of space, the constraint that the vector field $\nabla u$ does not vanish can be partially removed. Under appropriate conditions on the oscillations of the illumination $g$ on $\partial X$, stability results are obtained in \cite{A-AMPA-86} in cases where $\nabla u$ is allowed to vanish.

\subsection{Reconstructions from quadratic functionals in $u$}
\label{sec:TAT}

\subsubsection{Reconstructions under smallness conditions}
\label{sec:recsmall}

The TAT and (simplified) AOT problems are examples of a more general class we define as follows. Let $P(x,D)$ be an operator acting on functions defined in $\Cm^m$ for $m\in\Nm^*$ an integer and with values in the same space. Consider the equation
\begin{equation}
  \label{eq:systP}
  \begin{array}{rcll}
  P(x,D) u &=& \conductivitypot(x) u,\quad & x\in X \\
  u &=& f, \quad & x\in\partial X.
  \end{array}
\end{equation}
We assume that the above equation admits a unique weak solution in some Hilbert space $\mH_1$ for sufficiently smooth illuminations $f(x)$ on $\partial X$.

For instance, $P$ could be the Helmholtz operator $ik^{-1}(\Delta+k^2)$  seen in the preceding section with $u\in \mH_1:=H^1(X;\Cm)$ and $f\in H^\frac12(\partial X;\Cm)$.  Time-harmonic Maxwell's equations can be put in that framework with $m=n$ and
\begin{equation} \label{eq:maxP}
   P(x,D) = \dfrac{1}{ik}( \nabla\!\times\nabla\!\times  - k^2).
\end{equation}
We impose an additional constraint on $P(x,D)$ that the equation $P(x,D)u=f$ on $X$ with $u=0$ on $\partial X$ admits a unique solution in $\mH=L^2(X;\Cm^m)$. For instance, $\mH=L^2(X;\Cm)$ in the example seen in the preceding section in the scalar approximation provided that $k^2$ is not an eigenvalue of $-\Delta$ on $X$. For Maxwell's equations, the above constraint is satisfied so long as $k^2$ is not an internal eigenvalue of the Maxwell operator \cite{dlen3}. This is expressed by the existence of a constant $\alpha>0$ such that:
\begin{equation}
  \label{eq:constalpha}
  (P(x,D)u,u)_{\mH} \geq \alpha (u,u)_{\mH}.
\end{equation}
We assume that the conductivity $\conductivitypot$ is bounded from above by a positive constant:
\begin{equation}
  \label{eq:sigmaminmax}
   0< \conductivitypot(x) \leq \conductivitypot_M \qquad \mbox{ a.e. } x\in X.
\end{equation}
We denote by $\Sigma_{M}$ the space of functions $\conductivitypot(x)$ such that \eqref{eq:sigmaminmax} holds. Measurements are of the form $H(x)=\conductivitypot(x)|u|^2$, where $|\cdot|$ is the Euclidean norm on $\Cm^m$. Then we have the following result.
\begin{theorem}[\cite{BRUZ-IP-11}]
  \label{thm:smallsigma}
  Let $\conductivitypot_j\in \Sigma_{M}$ for $j=1,2$. Let $u_j$ be the solution to $P(x,D)u_j=\conductivitypot_j u_j$ in $X$ with $u_j=f$ on $\partial X$ for $j=1,2$. Define the internal functionals $H_j(x)=\conductivitypot_j(x)|u_j(x)|^2$ on $X$.

   Then for $\conductivitypot_M$ sufficiently small so that $\conductivitypot_M<\alpha$, we find that:
   \\
   (i) [Uniqueness] If $H_1=H_2$ $a.e.$ in $X$, then $\conductivitypot_1(x)=\conductivitypot_2(x)$ $a.e.$ in $X$ where $H_1=H_2>0$. \\
   (ii) [Stability] Moreover, we have the following stability estimate
   \begin{equation}
  \label{eq:stabestsmall}
  \|(\sqrt{\conductivitypot_1}-\sqrt{\conductivitypot_2})w_1\|_{\mH} \leq C \|(\sqrt{H_1}-\sqrt{H_2})w_2\|_{\mH},
\end{equation}
for some universal constant $C$ and for positive weights given by
\begin{equation}
  \label{eq:w12}
    w_1^2 (x) = \prod\limits_{j=1,2} \dfrac{|u_j|}{\sqrt\conductivitypot_j}(x) ,\quad w_2(x) =   \dfrac{1}{\alpha-\sup\limits_{x\in X}\sqrt{\conductivitypot_1\conductivitypot_2}}\max\limits_{j=1,2}\dfrac{\sqrt\conductivitypot_j}{|u_{j\rq{}}|}(x) + \max\limits_{j=1,2}\dfrac{1}{\sqrt\conductivitypot_j}(x) .
\end{equation}
Here $j'=j'(j)$ is defined as $j'(1)=2$ and $j'(2)=1$.
\end{theorem}
The theorem uses the spectral gap in \eqref{eq:constalpha}. Some straightforward algebra shows that
\begin{displaymath}
   P(x,D) (u_1-u_2) = \sqrt{\conductivitypot_1\conductivitypot_2} \big( |u_2|\hat u_1-|u_1|\hat u_2\big) + (\sqrt {H_1} - \sqrt{H_2}) \Big(\dfrac{\sqrt{\conductivitypot_1}}{|u_1|}-\dfrac{\sqrt{\conductivitypot_2}}{|u_2|}\Big).
\end{displaymath}
Here we have defined $\hat u = \frac{u}{|u|}$. Although this does not constitute an equation for $u_1-u_2$, it turns out that
\begin{displaymath}
    | |u_2|\hat u_1-|u_1|\hat u_2 | = | u_2-u_1|.
\end{displaymath}
This combined with \eqref{eq:constalpha} yields the theorem after some elementary manipulations \cite{BRUZ-IP-11}. 

\subsubsection{Reconstructions for the Helmholtz equation}
\label{sec:helm}

Let us consider the scalar model of TAT. We assume that $\conductivitypot\in H^p(X)$ for $p>\frac n2$ and construct
\begin{equation} \label{eq:q1}
  q(x)= k^2 + i k \conductivitypot(x) \in H^{p}(X), \qquad
  p>\frac n2.
\end{equation}
We assume that $q(x)$ is the restriction to $X$ of the compactly supported function (still called $q$) $q\in H^p(\Rm^n)$. The extension is chosen so that \cite{BU-IP-10}
\begin{math}
  \|q_{|X}\|_{H^{p}(X)} \leq C \|q\|_{H^{p}(\Rm^n)}
\end{math}
for some constant $C$ independent of $q$. Then \eqref{eq:invbrdy} may be recast as
\begin{equation}
  \label{eq:invbrdy2}
  \begin{array}{ll}
   \Delta u + q(x) u =0 \quad\mbox{ in } \quad X,\qquad  
   u=f \quad \mbox{ on } \quad \partial X.
  \end{array}
\end{equation}
The measurements are of the form $H(x)=\conductivitypot(x)|u|^2(x)$.

The inverse problem consists of reconstructing $q(x)$ from knowledge of $H(x)$. Note that $q(x)$ need not be of the form \eqref{eq:q1}. It could be a real-valued potential in a Helmholtz equation as considered  in \cite{T-IP-10} with applications in the so-called inverse medium problem. The reconstruction of $q(x)$ in \eqref{eq:invbrdy2} from knowledge of $H(x)=\conductivitypot(x)|u|^2(x)$ has been analyzed in \cite{BRUZ-IP-11,T-IP-10}. 

In the first reference, we have the following global stability reconstruction result.  We define $Y=H^p(X)$ for $p>\frac n2$ and $\mathcal M$ as the space of functions in $Y$
with norm bounded by a fixed (arbitrary) $M>0$. Let us define $Z=H^{p-\frac12}(\partial X)$. 
\begin{theorem}[\cite{BRUZ-IP-11}]
  \label{thm:reconstTAT}
  Let $\conductivitypot$ and $\tilde\conductivitypot$ be functions in $\mathcal M$.
  Let $f\in Z$ be a given (complex-valued) illumination and $H(x)$ be the measurement
  given in \eqref{eq:IScID} for $u$ solution of \eqref{eq:invbrdy}.
  Let $\tilde H(x)$ be the measurement constructed by replacing
  $\conductivitypot$ by $\tilde\conductivitypot$ in \eqref{eq:IScID} and
  \eqref{eq:invbrdy}.

  Then there is an open set of illuminations $f$ in $Z$ such that
  $H(x)=\tilde H(x)$ in $Y$ implies that $\conductivitypot(x)=\tilde\conductivitypot(x)$
  in $Y$. Moreover, there exists a constant $C$ independent of
  $\conductivitypot$ and $\tilde\conductivitypot$ in $\mathcal M$ such that
  \begin{equation}
    \label{eq:stab}
    \|\conductivitypot-\tilde\conductivitypot\|_Y\leq C \|H-\tilde H\|_Y.
  \end{equation}
\end{theorem}
The theorem is written in terms of $\conductivitypot$, which is the parameter of interest in TAT. The same result holds if $\conductivitypot$ is replaced by $q(x)$ in \eqref{eq:stab}. The reconstruction of $\conductivitypot$ is also constructive as the application of a Banach fixed point theorem. The proof is based on the construction of complex geometric optics solutions that will be presented in section \ref{sec:CGO}.

In reference \cite{T-IP-10}, the following local stability result is obtained. 
\begin{theorem}[\cite{T-IP-10}]
  \label{thm:triki}
  Let $q(x)\geq c_0>0$ be real-valued, positive, bounded on $X$ and such that $0$ is not an eigenvalue of $\Delta+q$ with domain $H^1_0(X)\cap H_2(\bar X)$. Let $\tilde q$ satisfy the same hypotheses and let $H$ and $\tilde H$ be the corresponding measurements. 
  
  Then there is a constant $\eps>0$ such that if $q$ and $\tilde q$ are $\eps-$close in $L^\infty(X)$ and if $f$ is in an $\eps-$dependent open set of (complex-valued) illuminations, then there is a constant $C$ such that 
 \begin{equation}
  \label{eq:stabtriki}
  \|q-\tilde q\|_{L^2(X)} \leq C \|H-\tilde H \|_{L^2(X)}.
\end{equation} 
\end{theorem}

Both theorems \ref{thm:reconstTAT} and \ref{thm:triki} show that the TAT and the inverse medium problem are stable inverse problems.  This is confirmed by the numerical reconstructions in \cite{BRUZ-IP-11}. The first result is more global but requires more regularity of the coefficients. It is based on the use of complex geometric optics CGO solutions to show that an appropriate functional is contracting in the space of continuous functions. The second result is more local in nature (a global uniqueness result is also proved in \cite{T-IP-10}) but requires less smoothness on the coefficient $q(x)$ and provides a stability estimate in the larger space $L^2(X)$. It also uses CGO solutions to show that the norm of a complex-valued solution to an elliptic equation is bounded from below by a positive constant. In both cases, the CGO solutions have traces at the boundary $\dX$ and the chosen illumination $f$ needs to be chosen in the vicinity of such traces. 

The results obtained in Theorem \ref{thm:smallsigma} under smallness constraints on $\conductivitypot$ apply for very general illuminations $f$. The above two results apply for more general (essentially arbitrary) coefficients but require more severe constraints on the illuminations $f$. 

\subsubsection{Non-unique reconstruction in the AOT setting}
\label{sec:nonuniq}

The above results concern the uniqueness of the reconstruction of the potential in a Helmholtz equation when well-chosen complex-valued boundary conditions are imposed. They also show that the reconstruction of $0<c_0\leq q(x)$  in $\Delta u+qu=0$ with real-valued $u=f$  from knowledge of $qu^2$ is unique. This corresponds to $P(x,D)=-\Delta$. In a simplified version of the acousto-optics problem considered in \cite{BS-PRL-10}, it is interesting to look at the problem where $P(x,D)=\Delta$ and where the measurements are given by $H(x)=\conductivitypot(x) u^2(x)$. Here, $u$ is the solution of the elliptic equation $(-\Delta + \conductivitypot) u =0$ on $X$ with $u=f$ on $\partial X$. Assuming that $f$ is non-negative, which is the physically interesting case, we obtain that $|u|=u$ and hence
\begin{displaymath}
  \Delta(u_1-u_2) = \sqrt{\conductivitypot_1\conductivitypot_2} (u_2-u_1) + (\sqrt{H_1}-\sqrt{H_2})\Big(\dfrac{\conductivitypot_1}{\sqrt{H_1}}-\dfrac{\conductivitypot_2}{\sqrt{H_2}}\Big).
 \end{displaymath}
 Therefore, as soon as $0$ is not an eigenvalue of $\Delta+\sqrt{\conductivitypot_1\conductivitypot_2} $, we obtain that $u_1=u_2$ and hence that $\conductivitypot_1=\conductivitypot_2$. For $\conductivitypot_0$ such that $0$ is not an eigenvalue of $\Delta+\conductivitypot_0$, we find that for $\conductivitypot_1$ and $\conductivitypot_2$ sufficiently close to $\conductivitypot_0$, then $H_1=H_2$ implies that $\conductivitypot_1=\conductivitypot_2$ on the support of $H_1=H_2$. 
 
However,  it is shown in \cite{BR-ProcAMS-11} that two different, positive, absorptions $\conductivitypot_{j}$ for $j=1,2$, may in some cases provide the same measurement $H=\conductivitypot_{j} u_{j}^2$ with $\Delta u_{j}=\conductivitypot_{j} u_{j}$ on $X$ with $u_{j}=f$ on $\partial X$ and in fact $\conductivitypot_1=\conductivitypot_2$ on $\partial X$ so that these absorptions cannot be  distinguished by their traces on $\partial X$. This counter-example shows that conditions such as the smallness condition in Theorem \ref{thm:smallsigma} are necessary in general. 
 
More generally, and following \cite{BR-ProcAMS-11}, consider an  elliptic problem of the form
\begin{equation}
  \label{eq:diff1}
     \begin{array}{ll}
     P u = \conductivitypot u \quad\mbox{ in }\quad X,\qquad   u=f \quad\mbox{ on }\quad \partial X,
 \end{array}
\end{equation}
and assume that  measurements of the form $H(x)=\conductivitypot(x)u^2(x)$ are available. Here, $P$ is a self-adjoint, non-positive, elliptic operator, which for concreteness we will take of the form $Pu=\nabla\cdot \diffusion(x)\nabla u$ with $\diffusion(x)$ known, sufficiently smooth, and bounded above and below by positive constants. We assume $f>0$ and $\diffusion>0$ so that by the maximum principle, $u>0$ on $X$. We also assume enough regularity on $\partial X$ and $f$ so that $u\in C^{2,\beta}(\bar X)$ for some $\beta>0$ \cite{gt1}. 

We observe that 
\begin{equation}
  \label{eq:nonlindiff}
     \begin{array}{ll}
     u P u = H \quad\mbox{ in }\quad X,\qquad   u=f \quad\mbox{ on }\quad \partial X
 \end{array}
\end{equation}
so that the inverse problem may be recast as a semilinear problem. The non-uniqueness result is an example of an Ambrosetti-Prodi result \cite{AP-AMPA-72} and in some sense generalizes the observation that $x\to x^2$ admits $0$, $1$, or $2$ (real-valued) solution(s) depending on the value of $x^2$. Let us define 
\begin{equation}
  \label{eq:phi}
  \phi: C^{2,\beta}(\bar X)\to C^{0,\beta}(\bar X), \qquad u \mapsto \phi(u) = u Pu.
\end{equation}
The singular points of $\phi$ are calculated from its first-order Fr\'echet derivative:
\begin{equation}
  \label{eq:Frechet}
  \phi'(u) v = v Pu + u Pv . 
\end{equation}
The operator $\phi'(u)$ is not invertible when $\conductivitypot:=\frac{Pu}{u}$ is such that $P+\lambda\conductivitypot$ admits $\lambda=1$ as an eigenvalue.  Let $\conductivitypot_0$ be such that $P+\conductivitypot_0$ is not invertible. We assume that the corresponding eigen-space is one dimensional and spanned by the eigenvector $\psi>0$ on $\dX$ such that $(P+\conductivitypot_0)\psi=0$  and $\psi=0$ on $\dX$. Let us define $u_0$ as
\begin{equation}
  \label{eq:u0}
  Pu_0 = \conductivitypot_0 u_0\,\mbox{ in }\, X,\qquad u_0=f \,\mbox{ on }\, \partial X,\qquad \conductivitypot_0>0.
\end{equation}
Moreover, $u_0$ is a singular point of $\phi(u)$ with $\phi'(u_0) \psi=0$ .
Then {\em define}
\begin{eqnarray}
  \label{eq:udelta}
  u_\delta& := &u_0+ \delta\psi, \quad X, \qquad \delta\in (-\delta_0,\delta_0)\\
  \label{eq:meas}
  \conductivitypot_\delta &:=& \dfrac{Pu_\delta}{u_\delta} = \conductivitypot_0 \dfrac{u_0-\delta\psi}{u_0+\delta\psi},\qquad
  H_\delta := \conductivitypot_\delta u_\delta^2 = \conductivitypot_0 u_\delta u_{-\delta} = \conductivitypot_0( u_0^2 -\delta^2 \psi^2).
\end{eqnarray}
We choose $\delta_0$ such that $\conductivitypot_\delta>0$ a.e. on $X$ for all $\delta\in  (-\delta_0,\delta_0)$.  Then, we have:
\begin{proposition}[\cite{BR-ProcAMS-11}]\label{prop:2}
Let $u_0$ be a singular point and $H_0=\phi(u_0)$ a critical value of $\phi$ as above and let 
$\psi$ be the normalized solution of $\phi'(u_0)\psi=0$. Let
$u_\delta$, $\conductivitypot_\delta$, and $H_\delta$ be defined as in \eqref{eq:udelta}-\eqref{eq:meas} for $0\not=\delta\in (-\delta_0,\delta_0)$ for $\delta_0$ sufficiently small.
Then we verify that:
\begin{displaymath} 
  \conductivitypot_\delta\not=\conductivitypot_{-\delta},\quad \conductivitypot_\delta>0, \quad H_\delta=H_{-\delta},\quad Pu_\delta=\conductivitypot_\delta u_\delta\,\,\mbox{ in } X,\quad u_\delta =f \,\, \mbox{ on } \partial X.
\end{displaymath}
\end{proposition}
This shows the non-uniqueness of the reconstruction of $\conductivitypot$ from knowledge of $H=\conductivitypot u^2$. Moreover we verify that $\conductivitypot_{\pm\delta}$ agree on $\partial X$ so that this boundary information cannot be used to distinguish between $\conductivitypot_\delta$ and $\conductivitypot_{-\delta}$. The non-uniqueness result is not very restrictive since we have seen that two coefficients, hence one coefficient, may be uniquely reconstructed from {\em two} well-chosen illuminations in the PAT results. Nonetheless, the above result shows once more that identifiability of the unknown coefficients is not always guaranteed by the availability of internal measurements.

\section{Reconstructions  from functionals of $\nabla u$}
\label{sec:diff}
  
We have seen two models of hybrid inverse problems with measurements involving $\nabla u$. In UMEIT, the measurements are of the form $H(x)=\conductivity(x)|\nabla u|^2(x)$ whereas in CDII, they are of the form $H(x)=\conductivity(x)|\nabla u|(x)$. 

Let us consider more generally measurements of the form $H(x)=\conductivity(x)|\nabla u|^{2-p}$ for $u$ the solution to the elliptic equation 
\begin{equation}
  \label{eq:ellipticsigma1}
  -\nabla\cdot \conductivity(x) \nabla u  = 0\quad \mbox{ in } X, \qquad u=f\quad \mbox{ on } \partial X.
\end{equation}
Since $H(x)$ is linear in $\conductivity(x)$, we have formally what appears to be an extension to $p\geq0$ of the $p-$Laplacian elliptic equations
\begin{equation}\label{eq:pLaplacian} 
   -\nabla\cdot \dfrac{H(x)}{|\nabla u|^{2-p}} \nabla u=0, 
\end{equation}
posed on a bounded, smooth, open domain $X\subset\Rm^n$, $n\geq2$, with prescribed Dirichlet conditions, say. When $1<p<\infty$, the above problem is known to admit a variational formulation with convex functional $J[\nabla u]=\int_X H(x) |\nabla u|^{p} dx$, which admits a unique minimizer in an appropriate functional setting solution of the above associated Euler-Lagrange equation \cite{evans}.

When $p=1$, the equation becomes degenerate while for $p<1$, the equation is in fact hyperbolic. We consider the problem of measurements that are quadratic (or bilinear) in $\nabla u$ (with applications to UMEIT and UMOT) in the next two sections. In the following section, we consider the case $p=1$.

%
%
\subsection{Reconstruction from a single power density measurement}
\label{sec:singpowerdensity}

The presentation follows \cite{B-UMEIT-11}. When $p=0$ so that measurements are of the form $H(x)=\conductivity(x)|\nabla u|^2$, the above $0-$Laplacian turns out to be a hyperbolic equation. Anticipating this behavior, we assume the availability of Cauchy data (i.e., $u$ and $\conductivity\nu\cdot\nabla u$ with $\nu$ the unit outward normal to $X$) on $\dX$ rather than simply Dirichlet data. Then \eqref{eq:pLaplacian} with $p=0$ becomes after some algebra
\begin{equation}
  \label{eq:Cauchy2}
  (I-2\widehat{\nabla u}\otimes\widehat{\nabla u}) : \nabla^2 u + \nabla \ln H\cdot\nabla u =0\,\mbox{ in } X,\qquad u=f \,\,\mbox{ and } \,\, \pdr{u}{\nu} = j \,\, \mbox{ on } \partial X.
\end{equation}
Here $\widehat{\nabla u}=\frac{\nabla u}{|\nabla u|}$. With 
\begin{equation}
  \label{eq:gh}
  g^{ij}=g^{ij}(\nabla u)=-\delta^{ij}+2(\widehat{\nabla u})_i(\widehat{\nabla u})_j\quad\mbox{ and } \quad
  k^i=-(\nabla \ln H)_i,
\end{equation}
the above equation is in coordinates 
\begin{equation}
  \label{eq:Cauchy}
    g^{ij}(\nabla u) \partial^2_{ij} u + k^i \partial_i u=0\,\mbox{ in } X,\qquad u=f \,\,\mbox{ and } \,\, \pdr{u}{\nu} = j \,\, \mbox{ on } \partial X.
\end{equation}
Since $g^{ij}$ is a definite matrix of signature $(1,n-1)$, then \eqref{eq:Cauchy} is a quasilinear strictly hyperbolic equation. The Cauchy data $f$ and $j$ then need to be provided on a space-like hyper-surface in order for the hyperbolic problem to be well-posed \cite{H-II-SP-83}. This is the main difficulty with solving \eqref{eq:Cauchy} with redundant Cauchy boundary conditions.

In general, we cannot hope to reconstruct $u(x)$, and hence $\conductivity(x)$ on the whole domain $X$. The reason is that the direction of ``time" in the second-order hyperbolic equation is $\widehat{\nabla u}$(x). The normal $\nu(x)$ at the boundary $\partial X$ will distinguish between the (good) part of $\partial X$ that is ``space-like" and the (bad) part of $\partial X$ that is ``time-like". Space-like surfaces such as $t=0$ provide stable information to solve the standard wave equation whereas in general it is known that arbitrary singularities can form in a wave equation from information on ``time-like" surfaces such as $x=0$ or $y=0$ in a three dimensional setting (where $(t,x,y)$ are local coordinates of $X$) \cite{H-II-SP-83}.

\subsubsection{Uniqueness and stability}
\label{sec:localuniq}

Let $(u,\conductivity)$ and $(\tilde u,\tilde\conductivity)$ be two solutions of the Cauchy problem \eqref{eq:Cauchy} with measurements $(H,f,j)$ and $(\tilde H,\tilde f,\tilde j)$, where we define the reconstructed conductivities
\begin{equation}
  \label{eq:reccond} \conductivity(x) = \dfrac{H}{|\nabla u|^2}(x), \qquad \tilde \conductivity(x) = \dfrac{\tilde H}{|\nabla \tilde u|^2}(x).
\end{equation}
Let $v=\tilde u- u$. We find that 
\begin{displaymath} 
 \nabla\cdot \Big( \dfrac{H}{|\nabla \tilde u|^2|\nabla u|^2} \Big\{ (\nabla u+\nabla \tilde u) \otimes (\nabla u+\nabla \tilde u)-(|\nabla u|^2+|\nabla \tilde u|^2)I \Big\} \nabla v + \delta H \Big(\dfrac{\nabla \tilde u}{|\nabla \tilde u|^2}+\dfrac{\nabla u}{|\nabla u|^2}\Big)\Big)=0.
\end{displaymath}
This equation is recast as
\begin{equation}
  \label{eq:linCauchy} 
  \mg^{ij}(x) \partial^2_{ij} v + \mk^i \partial_i v + \partial_i (l^i \delta H)=0 \quad \mbox{ in } X,\qquad
  v= \tilde f-f, \quad \pdr{v}{\nu} = \tilde j-j \quad \mbox{ on } \partial X,
\end{equation}
for appropriate coefficients  $\mk^i$ and $l^i$, where
\begin{equation}\label{eq:mg} 
  \begin{array}{rcl}
  \mg(x) &=&  \dfrac{H}{|\nabla \tilde u|^2|\nabla u|^2} \Big\{ (\nabla u+\nabla \tilde u) \otimes (\nabla u+\nabla \tilde u)-(|\nabla u|^2+|\nabla \tilde u|^2)I \Big\} \\[3mm]
    &=& 
    \alpha(x) \Big(\,\be(x) \otimes \be(x) - \beta^2(x) \big(I-\be(x) \otimes \be(x) \big) \Big),
    \end{array}
 \end{equation} 
with
\begin{equation}\label{eq:bebeta}
    \be(x) = \dfrac{\nabla u+\nabla \tilde u}{|\nabla u+\nabla \tilde u|}(x) ,\qquad \beta^2(x) = \dfrac{|\nabla u+\nabla \tilde u|^2}{|\nabla u+\nabla \tilde u|^2- (|\nabla u|^2+|\nabla \tilde u|^2)}(x) ,
 \end{equation}
and $\alpha(x)$ is the appropriate (scalar) normalization coefficient. For $\nabla u$ and $\nabla\tilde u$ sufficiently close so that $\nabla u\cdot\nabla \tilde u>0$, then the above linear equation for $v$ is strictly hyperbolic. We define the Lorentzian metric $\mh=\mg^{-1}$ so that $\mh_{ij}$ are the coordinates of the inverse of the matrix $\mg^{ij}$. We denote by $\aver{\cdot,\cdot}$ the bilinear product associated to $\mh$ so that $\aver{u,v}=\mh_{ij}u^iv^j$ where the two vectors $u$ and $v$ have coordinates $u^i$ and $v^i$, respectively. We verify that 
\begin{equation}
  \label{eq:mh}
  \mh(x) \,=\, \frac1{\alpha(x)} \Big(\,\be(x) \otimes \be(x) - \dfrac1{\beta^{2}(x)} \big(I-\be(x) \otimes \be(x) \big) \Big).
\end{equation}
 
 The space-like part $\Sigma_g$ of $\partial X$ is given by $\mh(\nu,\nu)>0$, i.e., $\nu$ is a time-like vector, or equivalently
\begin{equation}
  \label{eq:spacelikepartialX}
  |\nu(x) \cdot \be(x) |^2 >  \frac 1{1+\beta^2(x)} \qquad x\in \partial X.
\end{equation}
Above, the ``dot" product is with respect to the standard Euclidean metric and $\nu$ is a unit vector for the Euclidean metric, not for the metric $\mh$. Let $\Sigma_1$ be an open connected component of $\Sigma_g$ and let $\cO=\cup_{0<\tau<s}\Sigma_2(\tau)$ be a domain of influence of $\Sigma_1$ swept out by the space-like surfaces $\Sigma_2(\tau)$; see \cite{B-UMEIT-11,Taylor-PDE-1}. Then we have the following local stability result:
\begin{theorem}[Local Uniqueness and Stability.]
  \label{thm:LinearStab}
  Let $u$ and $\tilde u$ be two solutions of \eqref{eq:Cauchy}.
  We assume that $\mg$ constructed in \eqref{eq:mg} is strictly hyperbolic. Let $\Sigma_1$ be an open connected component of $\Sigma_g$ the space-like component of $\partial X$ and let $\cO$ be a domain of influence of $\Sigma_1$ constructed as above. Let us define the energy
  \begin{equation}
  \label{eq:energy}
  E(dv) = \aver{dv,\nu_2}^2 - \dfrac12 \aver{dv,dv}\aver{\nu_2,\nu_2}.
\end{equation}
Here, $dv$ is the gradient of $v$ in the metric $\mh$ given in coordinates by $\mg^{ij}\partial_j v$. Then:
\begin{equation}
  \label{eq:localstab}
  \dint_{\cO} E(dv) dx \leq C \Big( \dint_{\Sigma_1} |f-\tilde f|^2 + |j-\tilde j|^2 \,d\sigma + \dint_{\cO} |\nabla \delta H|^2 \,dx\Big),
\end{equation}
where $dx$ and $d\sigma$ are the standard measures on $\cO$ and $\Sigma_1$, respectively. 

In the Euclidean metric, let $\nu_2(x)$ be the unit vector to $x\in\Sigma_2(\tau)$, define $c(x) := \nu_2(x) \cdot \be(x)$ and
\begin{equation}
  \label{eq:theta} \theta:= \min_{x\in \Sigma_2(\tau)} \Big[c^2(x) - \dfrac{1}{1+\beta^2(x)}\Big].
\end{equation}
Then we have that 
\begin{equation}
  \label{eq:localstab2}
  \dint_{\cO} |v^2|+|\nabla v|^2 + (\conductivity-\tilde\conductivity)^2 \,dx \leq \dfrac{C}{\theta^2} \Big( \dint_{\Sigma_1} |f-\tilde f|^2 + |j-\tilde j|^2 \,d\conductivity + \dint_{\cO} |\nabla \delta H|^2 \,dx\Big),
\end{equation}
where $\conductivity$ and $\tilde\conductivity$ are the conductivities in \eqref{eq:reccond}.
Provided that $f=\tilde f$, $j=\tilde j$, and $H=\tilde H$, we obtain that $v=0$ and the uniqueness result $u=\tilde u$ and $\conductivity=\tilde\conductivity$.
\end{theorem}
The proof is based on adapting energy methods for hyperbolic equations as they are summarized in \cite{Taylor-PDE-1}. The energy $E(dv)$ fails to control $dv$ for null-like or space-like vectors, i.e., $\mh(dv,dv)\leq0$. The parameter $\theta$ measures how time-like the vector $dv$ is on the domain of influence $\cO$.  As $\cO$ approaches the boundary of the domain of influence of $\Sigma_g$ and $\theta$ tends to $0$, the energy estimates deteriorate as indicated in \eqref{eq:localstab2}.

Assuming that the errors on the Cauchy data $f$ and $j$ are negligible, we obtain the following stability estimate for the conductivity
\begin{equation}
  \label{eq:stabcond1}
    \|\conductivity-\tilde\conductivity\|_{L^2(\cO)} \leq \dfrac{C}{\theta} \|H-\tilde H\|_{H^1(X)}.
\end{equation}
Under additional regularity assumptions on $\conductivity$, for instance assuming that $H\in H^s(X)$ for $s\geq2$, we find by standard interpolation that
\begin{equation}
  \label{eq:stabcond2}
    \|\conductivity-\tilde\conductivity\|_{L^2(\cO)} \leq \dfrac{C}{\theta} \|H-\tilde H\|^{1-\frac 1s}_{L^2(X)} \|H+\tilde H\|^{\frac1 s}_{H^s(X)},
\end{equation}
We thus obtain H\"older-stable reconstructions in the practical setting of square integrable measurement errors. However, stability is {\em local}. Only on the domain of influence of the space-like part of the boundary can we obtain a stable reconstruction. This can be done by solving a nonlinear strictly hyperbolic equation analyzed in \cite{B-UMEIT-11} using techniques summarized in \cite{H-SP-97}. 

\subsubsection{Global reconstructions}
\label{sec:global0lap}

In the preceding result, the main roadblock to global reconstructions was that the domain of influence of the space-like part of the boundary was a strict subset of $X$. There is a simple solution to this problem: simply make sure that the whole boundary is a level set of $u$ and that no critical points of $u$ (where $\nabla u=0$) exist. Then all of $X$ is in the domain of influence of the space-like part of $\dX$, which is the whole of $\dX$. This setting can be made possible independent of the conductivity $\conductivity$  in two dimensions of space but not always in higher dimensions. 

Let $n=2$. We assume that $X$ is an open smooth domain diffeomorphic to an annulus with boundary $\partial X=\partial X_0\cup\partial X_1$. We assume that $f=0$ on the external boundary $\partial X_0$ and that $f=1$ on the internal boundary $\partial X_1$. The boundary of $X$ is thus composed of two smooth connected components that are different level sets of the solution $u$. The solution $u$ to \eqref{eq:ellipticsigma1} is uniquely defined on $X$. Then we can show:
\begin{proposition}[\cite{B-UMEIT-11}]
\label{prop:doughnut2d}
We assume that both the geometry of $X$ and $\conductivity(x)$ are sufficiently smooth. Then $|\nabla u|$ is bounded from above and below by positive constants. The level sets $\Sigma_c=\{x\in X,\, u(x)=c\}$ for $0<c<1$ are smooth curves that separate $X$ into two disjoint subdomains.
\end{proposition}
The proof is based on the fact that critical points of solutions to elliptic equations in two dimensions are isolated \cite{A-AMPA-86}. The result extends to higher dimensions provided that $|\nabla u|$ does not vanish with exactly the same proof. In the absence of critical points, we thus obtain that $\be(x)=\widehat{\nabla u}=\nu(x)$ so that $\nu(x)$ is clearly a time-like vector. Then the local results of Theorem \ref{thm:LinearStab} become global results, which yields the following proposition:
\begin{proposition}
\label{prop:doughnutnd}
 Let $X$ be the geometry described above in dimension $n\geq2$ and $u(x)$ the solution to \eqref{eq:ellipticsigma1}. We assume here that both the geometry and $\conductivity(x)$ are sufficiently smooth. We also {\em assume} that  $|\nabla u|$ is bounded from above and below by positive constants. Then the nonlinear equation \eqref{eq:Cauchy} admits a unique solution and the reconstruction of $u$ and of $\conductivity$ is stable in $X$ in the sense  described in Theorem \ref{thm:LinearStab}.
\end{proposition}

In dimensions $n\geq3$, we cannot guaranty that $u$ does not have any critical point independent of the conductivity. If the conductivity is close to a constant, then by continuity, $u$ does not have any critical point and the above result applies. This proves the result for sufficiently small perturbations of the case $\conductivity(x)=\conductivity_0$. In the general case, however, we cannot guaranty that $\nabla u$ does not vanish and in fact can produce counter-examples (see \cite{B-UMEIT-11}):
\begin{proposition}[\cite{B-UMEIT-11,BMN-ARMA-04,M-PAMS-93}]
\label{prop:badsigma3d}
   There is an example of a smooth conductivity such that $u$ admits critical points.
\end{proposition}
So in dimensions $n\geq3$, we are not guaranteed that the nonlinear equation will remain strictly hyperbolic. What we can do, however, is again to use the notion of complex geometric optics solutions. We have the result:
\begin{theorem}[\cite{B-UMEIT-11}]
  \label{thm:globalCGO}
  Let $\conductivity$ be extended by $\conductivity_0=1$ on $\Rm^n\backslash \tilde X$, where $\tilde X$ is the domain where $\conductivity$ is not known. We assume that $\conductivity$ is smooth on $\Rm^n$. Let $\conductivity(x)-1$ be supported without loss of generality on the cube $(0,1)\times(-\frac12,\frac12)^{n-1}$. Define the domain $X=(0,1)\times B_{n-1}(a)$, where $B_{n-1}(a)$ is the $n-1$-dimensional ball of radius $a$ centered at $0$ and where $a$ is sufficiently large that the light cone for the Euclidean metric emerging from $B_{n-1}(a)$ strictly includes $\tilde X$. Then there is an open set of illuminations $(f_1,f_2)$ such that if $u_1$ and $u_2$ are the corresponding solutions of \eqref{eq:ellipticsigma1}, then the following measurements
  \begin{equation}
  \label{eq:threemeas} H_{1}(x)=\conductivity(x) |\nabla u_1|^2(x),\quad H_{2}(x)=\conductivity(x)|\nabla u_2|^2(x),\quad
  H_{3}(x)= \conductivity(x)|\nabla (u_1+u_2)|^2,
\end{equation}
with the corresponding Cauchy data $(f_1,j_1)$, $(f_2,j_2)$ and $(f_1+f_2,j_1+j_2)$ at $x_1=0$ uniquely determine $\conductivity(x)$. Moreover, let $\tilde H_{i}$ be measurements corresponding to $\tilde\conductivity$ and $(\tilde f_1,\tilde j_1)$ and $(\tilde f_2,\tilde j_2)$ the corresponding Cauchy data at $x_1=0$. We assume that $\conductivity(x)-1$ and $\tilde\conductivity(x)-1$ are smooth and such that their norm in $H^{\frac n2+3+\eps}(\Rm^n)$ for some $\eps>0$ are bounded by $M$. Then for a constant $C$ that depends on $M$, we have the global stability result
\begin{equation}
  \label{eq:globalstab}
  \|\conductivity-\tilde\conductivity\|_{L^2(\tilde X)}\leq C \Big(\|d_C-\tilde d_C\|_{(L^2(B_{n-1}(a)))^4} + \sum_{i=1}^3 \|\nabla H_{i}-\nabla \tilde H_{i}\|_{L^2(X)}\Big).
\end{equation}
Here, we have defined 
$d_C=(f_1,j_1,f_2,j_2)$ with $\tilde d_C$ being defined similarly.
\end{theorem}
The ``three'' measurements $H_{i}$ in \eqref{eq:threemeas} actually correspond to ``two'' physical measurements since $H_{3}$  may be determined from the experiments yielding $H_{1}$ and $H_{2}$ as demonstrated in \cite{B-UMEIT-11,KK-AET-11}. The three measurements are constructed so that two independent strictly hyperbolic Lorentzian metrics can be construct everywhere inside the domain. These metrics are constructed by means of CGO solutions. The boundary conditions $f_j$ have to be close to the traces of the CGO solutions. We thus obtain a global Lipschitz stability result. The price to pay is that the open set of illuminations is not very explicit and may depend on the conductivities one seeks to reconstruct. 

\medskip

For conductivities that are close to a constant, several reconstructions are therefore available. We have seen that geometries of the form of an annulus (with a hole that can be arbitrarily small and arbitrarily close to the boundary where $f=0$) allowed us to obtain globally stable reconstructions since in such situations, it is relatively easy to avoid the presence of critical points. The method of CGO solutions can be shown to apply for a well-defined set of illuminations since the (harmonic) CGO solutions are explicitly known for the Euclidean metric and of the form $e^{\rho\cdot x}$ for $\rho$ a complex valued vector such that $\rho\cdot\rho=0$. After linearization in the vicinity of the Euclidean metric, another explicit reconstruction procedure was introduced in \cite{KK-AET-11}.

%
%
\subsection{Reconstructions from multiple power density measurements}
\label{sec:multpowerdensity}

Rather than reconstructing $\gamma$ from one given measurement of the form $\gamma(x)|\nabla u|^2$, we can instead acquire several measurements of the form
\begin{equation}
  \label{eq:bilinear0Lap}
  H_{ij}(x) = \conductivity(x) \nabla u_i (x)\cdot\nabla u_j(x) \quad \mbox{ in } X,\qquad 1\leq i,j\leq M,
\end{equation}
where $u_j$ solves the elliptic problem \eqref{eq:ellipticsigma1} with $f$ given by $f_j$ for $1\leq j\leq M$. The result presented in Theorem \ref{thm:globalCGO} above provides a positive answer for $M=2$ when the available internal functionals are augmented by Cauchy data at the boundary of the domain of interest.

Results obtained in \cite{BBMT-11,CFGK-SJIS-09,MB-IPI-12} and based on an entirely different procedure and not requiring knowledge of boundary data show that $M=2\lfloor \frac{n+1}2\rfloor$ measurements allow for a global reconstruction of $\conductivity$, i.e., $M=n$ for $n$ even and $M=n+1$ for $n$ odd. Such results were first obtained in \cite{CFGK-SJIS-09} in the case $n=2$ and have been extended with a slightly different presentation to the cases $n=2$ and $n=3$ in \cite{BBMT-11} while the general case $n\geq2$ is treated in \cite{MB-IPI-12}.
Let us assume that $n=3$ for concreteness. Then $H_{ij}=S_i\cdot S_j$, where we have defined 
\begin{displaymath} 
   S_j(x) = \sqrt{\conductivity(x)} \nabla u_j(x), \qquad 1\leq j\leq M.
\end{displaymath}
Let $S=(S_1,S_2,S_3)$ be a matrix of $n=3$ column vectors $S_j$.  Then $S^TS=H$, where $S^T$ is the transpose matrix made of the rows given by the $S_j$. We do not know $S$ or the $S_j$, but we know its normal matrix $S^TS=H$.  Let $T$ be a matrix such that $R=ST^T$ is a rotation-valued field on $X$. Two examples are $T=H^{-\frac12}$ or the lower-triangular $T$ obtained by the Gram-Schmidt procedure. We thus have information on $S$. We need additional equations to solve for $S$, or equivalently $R$, uniquely. The elliptic equation may be written as $\nabla\cdot \sqrt\conductivity S_j=0$, or equivalently
\begin{equation}\label{SjF} 
   \nabla\cdot S_j + F \cdot S_j=0,\qquad F=\nabla (\log\sqrt\conductivity) = \frac12\nabla \log\conductivity.
\end{equation}
Now, since $\conductivity^{-\frac{1}2}S_j$ is a gradient, its curl vanishes and we find that 
\begin{equation}
  \label{eq:curlSj}
  \nabla\times S_j - F\times S_j=0.
\end{equation}
Here, $F$ is unknown. We first eliminate it from the equations and then find a closed form equation for $S$ or equivalently for $R$ as a field in $SO(n;\Rm)$.

Let $T$ be the aforementioned matrix $T$, say $T=H^{-\frac12}$ with entries $t_{ij}$ for $1\leq i,j\leq n$. Let $t^{ij}$ be the entries of $T^{-1}$ and define the vector fields
\begin{equation}
  \label{eq:Vijl}
  V_{ij} := \nabla (t_{ik}) t^{kj}, \qquad \mbox{ i.e., } \quad V_{ij}^l = \partial_l(t_{ij})t^{kl},\quad 1\leq i,j,l\leq n.
\end{equation}
We then define $R(x)=S(x) T^T(x)\in SO(n;\Rm)$ the matrix whose columns are composed of the column vectors $R_j=S_j T^T$. Then in all dimension $n\geq2$, we find
\begin{lemma}[\cite{BBMT-11,MB-IPI-12}]\label{lem:F} In $n\geq2$, we have the following expression:
\begin{equation}
  \label{eq:FVij} 
  F=\dfrac1n \Big( \dfrac12 \nabla\log\det H + \dsum_{i,j=1}^n \big((V_{ij}+V_{ji})\cdot R_i\big)R_j\Big).
\end{equation}
\end{lemma}
The proof in dimension $n=2,3$ can be found in\cite{BBMT-11} and in arbitrary dimension in \cite{MB-IPI-12}.

Note that the determinant of $H$ needs to be positive on the domain $X$ in order for the above expression for $F$ to make sense. It is, however, difficult to ensure that the determinant of several gradients remains positive and there are in fact counter-examples as shown in \cite{BMN-ARMA-04}. Here again, complex geometric optics solutions are useful to control the determinant of gradients of elliptic solutions locally and globally using several solutions. We state a global result in the practical setting $n=3$.

Let be $m\ge 3$ solutions of the elliptic equation and assume that there exists an open covering ${\mathcal O} = \{\Omega_k\}_{1\le k\le N}$ ($X\subset\cup_{k=1}^N \Omega_i$), a constant $c_0>0$ and a function $\tau:[1,N]\ni i\mapsto \tau (i) = (\tau(i)_1, \tau(i)_2, \tau(i)_3) \in[1,m]^3$, such that
\begin{align}
    \inf_{x\in\Omega_i} \det(S_{\tau(i)_1}(x), S_{\tau(i)_2}(x), S_{\tau(i)_3}(x)) \ge c_0, \quad 1\le i\le N.
    \label{eq:condji}
\end{align}
Then we have the following result:
\begin{theorem}[3D global uniqueness and stability] \label{thm:stab3d}
    Let $X\subset\Rm^3$ be an open convex bounded set, and let two sets of $m\ge 3$ solutions of \eqref{eq:ellipticsigma1} generate measurements $(H, \tilde H)$ whose components belong to $W^{1,\infty}(X)$. Assume that one can define a couple $({\mathcal O}, \tau)$ such that \eqref{eq:condji} is satisfied for both sets of solutions $S$ and $\tilde S$. Let also $x_0\in\overline{\Omega}_{i_0}\subset\overline{X}$ and $\conductivity(x_0)$, $\tilde \conductivity(x_0)$, $\{S_{\tau(i_0)_i}(x_0), \tilde S_{\tau(i_0)_i}(x_0)\}_{1\le i\le 3}$ be given. 
Let $\conductivity$ and $\tilde\conductivity$ be the conductivities corresponding to the measurements $H$ and $\tilde H$, respectively. Then we have the following stability estimate:
    \begin{align}
	\|\log\conductivity-\log\tilde\conductivity\|_{W^{1,\infty}(X)} \le C\big( \epsilon_0 + \|H-\tilde H\|_{W^{1,\infty}(X)}\big),
	\label{eq:globstab3d}
    \end{align}
    where $\epsilon_0$ is the error at the initial point $x_0$
    \begin{align*}
	\epsilon_0 = |\log\conductivity_0-\log\tilde\conductivity_0| + \sum_{i=1}^3 \|S_{\tau(i_0)_i}(x_0) - \tilde S_{\tau(i_0)_i}(x_0) \|.
    \end{align*} 
\end{theorem}
This shows that the reconstruction of $\conductivity$ is stable from such redundant measurements. Moreover, the reconstruction is constructive. Indeed, after eliminating $F$ from the equations for $R$, we find an equation of the form $\nabla R = G(x,R)$, where $G(x, R)$ is polynomial of degree three in the entries of $R$. This is a redundant equation whose solution, when it exists, is unique and stable with respect to perturbations in $G$ and the conditions at a given point $x_0$.

That \eqref{eq:condji} is satisfied can again be proved by means of complex geometric optics solutions as is briefly mentioned in section \ref{sec:CGO} below; see \cite{BBMT-11}.

%
%
\subsection{Reconstruction from a single current density measurement}
\label{sec:singcurrentdensity}

Let us now come back to the $1-$Laplacian, which is a degenerate elliptic problem. In many cases, this problem admits multiple admissible solutions \cite{KKSY-SIMA-02}. The inverse problem then cannot be solved uniquely. In some settings, however, uniqueness can be restored \cite{KKSY-SIMA-02,NTT-IP-07,NTT-IP-09,NTT-Rev-11}. 

Recall that the measurements are of the form $H(x)=\conductivity |\nabla u|$ so that $u$ solves the following degenerate quasi-linear equation
\begin{equation}
  \label{eq:1laplace}
  \nabla\cdot \dfrac{H(x)}{|\nabla u|}\nabla u=0 \quad \mbox{ in } X.
\end{equation}
Different boundary conditions may then be considered. It is shown in \cite{KKSY-SIMA-02} that the above equation augmented with Neumann boundary conditions of the form
\begin{displaymath} 
    \dfrac{H}{|\nabla u|} \pdr u\nu =h \mbox{ on } \dX,\quad \dint_{\dX} u d\sigma =0,
\end{displaymath}
admits an infinite number of solutions once it admits a solution, and may also admit no solution at all. One possible strategy is to acquire two measurements of the form $H(x)=\conductivity |\nabla u|$ corresponding to two prescribed currents. In this setting, it is shown in \cite{KKSY-SIMA-02} that (appropriately defined) singularities of $\conductivity$ are uniquely determined by the measurements. We refer the reader to the latter reference for the details.

Alternatively, we may augment the above equation \eqref{eq:1laplace} with Dirichlet data. Then the reconstruction of $\conductivity$ was shown to be uniquely determined in \cite{NTT-IP-07,NTT-IP-09,NTT-Rev-11}.  Why Dirichlet conditions help to stabilize the equation may be explained as follows. The $1-$Laplace equation \eqref{eq:1laplace} may be recast as 
\begin{displaymath} 
   (I-\widehat{\nabla u}\otimes\widehat{\nabla u}):\nabla^2 u + \nabla\ln H\cdot\nabla u=0,
\end{displaymath}
following similar calculations to those leading to \eqref{eq:Cauchy2}. The only difference is the ``$2$'' in front of $\widehat{\nabla u}\otimes\widehat{\nabla u}$ replaced by ``$1$'', or more generally $2-p$ for a $p-$Laplacian. When $p>1$, the problem remains strictly elliptic. When $p<1$, the problem is hyperbolic, and when $p=1$, it is degenerate in the direction $\widehat{\nabla u}$ and elliptic in the transverse directions. We can therefore modify $u$ so that its level sets remain unchanged and still satisfy the above partial differential equation. This modification can also be performed so that Neumann boundary conditions are not changed. This is the procedure used in \cite{KKSY-SIMA-02} to show the non-uniqueness of the reconstruction for the 1-Laplacian with Neumann boundary conditions. 

Dirichlet conditions, however, are modified by changes in the level sets of $u$. It turns out that even with Dirichlet conditions, several (viscosity) solutions to \eqref{eq:1laplace} may be constructed when $H\equiv1$; see, e.g., \cite{NTT-IP-07,NTT-Rev-11}. However, such solutions involve vanishing gradients on sets of positive measure. 

The right formulation for the CDII inverse problem that allows one to avoid vanishing gradients is to recast \eqref{eq:1laplace} as the minimization of the functional
\begin{equation}
  \label{eq:funct1lap}
  F[\nabla v] = \dint_X H(x) |\nabla v| dx,
\end{equation}
over $v\in H^1(X)$ with $v=f$ on $\dX$. Note that $F[\nabla v]$ is convex although it is not {\em strictly} convex. Moreover, let $\conductivity$ be the conductivity and $H=\conductivity|\nabla u|$ the corresponding measurement. Let then $v\in H^1(X)$ with $v=f$ on $\dX$. Then,
\begin{displaymath} 
   F[\nabla v] = \dint_X \conductivity |\nabla u||\nabla v| dx \geq \dint_X \conductivity \nabla u \cdot \nabla v dx = \dint_{\dX} \sigma \pdr u\nu f ds = F[\nabla u],
\end{displaymath}
by standard integrations by parts. This shows that $u$ minimizes $F$. We have the following result:
\begin{theorem}[\cite{NTT-Rev-11}]
  \label{thm:cdii} Let $(f,H)\in C^{1,\alpha}(\dX)\times C^\alpha(\bar X)$ with $H=\conductivity|\nabla u|$ for some $\conductivity\in C^\alpha(\bar X)$. Assume that $H(x)>0$ a.e. in $X$. Then the minimization of 
  \begin{equation}
  \label{eq:mincdii} \mbox{ argmin } \big\{ F[\nabla v] ,\,\, v\in W^{1,1}(X)\cap C(\bar X),\\,\,v_{|\dX}=f \big\},
\end{equation}
has a unique solution $u_0$. Moreover $\sigma_0=H|\nabla u_0|^{-1}$ is the unique conductivity associated to the measurement $H(x)$. 
\end{theorem}

It is known in two dimensions of space that $H(x)>0$ is satisfied for a large class of boundary conditions $f(x)$; see Lemma \ref{lem:gradpos} in the next section. In three dimensions of space, however, critical points of $u$ may arise as observed earlier in this paper; see e.g., \cite{B-UMEIT-11}. The CGO solutions that are analyzed in the following section allow us to show that $H(x)>0$ holds for an open set of illuminations $f$ at the boundary of the domain $\partial X$; see \eqref{eq:lowernablau1} below. 


Several reconstruction algorithms have been devised in \cite{KKSY-SIMA-02,NTT-IP-07,NTT-IP-09,NTT-Rev-11}, to which we refer for additional details. The numerical simulations presented in these papers show that when uniqueness is guaranteed, then the reconstructions are very high resolution and quite robust with respect to noise in the data, as is expected for general hybrid inverse problems. 

%
\section{Qualitative properties of forward solutions}
\label{sec:propsol}

%
\subsection{The case of two spatial dimensions}
\label{sec:2d}

Several explicit reconstructions obtained in hybrid inverse problems require that the solutions to the considered elliptic equations satisfy specific qualitative properties such as the absence of any critical point  or the positivity of the determinant of gradients of solutions. Such results can be proved in great generality in dimension $n=2$ but do not always hold in dimension $n\geq3$. 

In dimension $n=2$, the critical points of $u$ (points $x$ where $\nabla u(x)=0$) are necessarily isolated as is shown in, e.g., \cite{A-AMPA-86}. From this and techniques of quasiconformal mappings that are also restricted to two dimensions of space, we can show the following results.

\begin{lemma}[\cite{AN-ARMA-01}]
Let $u_1$ and $u_2$ be the solutions of \eqref{eq:ellipticsigma1} on $X$ simply connected with boundary conditions $f_1=x_1$ and $f_2=x_2$ on $\dX$, respectively, where $x=(x_1,x_2)$ are Cartesian coordinates on $X$. Assume that $\conductivity$ is sufficiently smooth. Then $(x_1,x_2)\mapsto(u_1,u_2)$ from $X$ to its image is a diffeomorphism. In other words, $\det(\nabla u_1,\nabla u_2)>0$ uniformly on $\bar X$.
\end{lemma}
This result is useful in the analysis of UMEIT and UMOT in the case of redundant measurements. It is shown in \cite{BMN-ARMA-04} that the appropriate extension of this result is false in dimension $n\geq3$.

We recall that a function continuous on a simple closed contour is almost two-to-one if it is two-to-one except possibly at its maximum and minimum \cite{NTT-IP-07}. Then we have, quite similarly to the result in Lemma \ref{lem:crit} and Proposition \ref{prop:doughnut2d}, which also use the results in \cite{A-AMPA-86}, the following:
\begin{lemma}[\cite{NTT-IP-07}]\label{lem:gradpos}
  Let $X$ be a simply connected planar domain and let $u$ be solution of \eqref{eq:ellipticsigma1} with $f$ almost two-to-one and $\sigma$ sufficiently smooth. Then $|\nabla u|$ is bounded from below by a positive constant on $\bar X$. Moreover, the level sets of $u$ are open curves inside $X$ with their two end points on $\dX$.
\end{lemma}
This shows that for a large class of boundary conditions with one maximum and one minimum, the solution $u$ cannot have any critical point in $\bar X$. On an annulus with boundaries equal to level sets of $u$, we saw in Proposition \ref{prop:doughnut2d} that $u$ had no critical points on $X$ in dimension $n=2$. This was used to show that the normal vector to the level sets of $u$ always forms a time-like vector for the Lorentzian metric defined in \eqref{eq:Cauchy}.


\medskip

All these results no longer hold in dimension $n\geq3$. See, e.g., \cite{B-UMEIT-11,BMN-ARMA-04} for counter-examples. In dimension $n\geq3$, the required qualitative properties cannot be obtained for a given set of illuminations (boundary conditions) independent of the conductivity. However, for conductivities that are bounded (with an arbitrary bound) in an appropriate norm, there are open sets of illuminations that allow us to obtain the required qualitative properties. One way to construct such solutions is by means of the complex geometric optics solutions that are analyzed in the next section.

%
\subsection{Complex Geometric Optics solutions}
\label{sec:CGO}

\subsubsection{CGO solutions and Helmholtz equations.}
Complex geometrical optics (CGO) solutions allow us to treat the potential $q$ in the equation 
\begin{equation}\label{eq:HelCGO} 
    (\Delta +q)u=0 \mbox{ in } X,\qquad u=f \mbox{ on } \dX,
 \end{equation}
as a perturbation of the leading operator $\Delta$. When $q=0$, CGO solutions are harmonic solutions defined on $\Rm^n$ and are of the form
\begin{displaymath} 
     u_\rho(x) = e^{\rho\cdot x},\qquad \rho\in\Cm^n \mbox{ such that } \rho\cdot\rho=0.
 \end{displaymath}
For $\rho=\rho_r+i\rho_i$ with  $\rho_r$ and $\rho_i$ vectors in $\Rm^n$, this means that $|\rho_r|^2=|\rho_i|^2$ and $\rho_r\cdot\rho_i=0$.

When $q\not\equiv0$, CGO solutions are solutions of the following problem
\begin{equation}
  \label{eq:CGO0}
  \Delta u_\rho + qu_\rho =0 ,\qquad u_\rho \sim e^{\rho\cdot x} \mbox{ as } 
    |x|\to\infty.
\end{equation}

More precisely, we say that $u_\rho$ is a solution of the above equation
with $\rho\cdot\rho=0$ and the proper behavior at infinity when it is
written as
\begin{equation}
  \label{eq:decurho}
  u_\rho(x) = e^{\rho\cdot x} \big(1+\psi_\rho(x)\big),
\end{equation}
for $\psi_\rho\in L^2_\delta$ a weak solution of
\begin{equation}
  \label{eq:psi}
  \Delta \psi_\rho + 2\rho\cdot\nabla \psi_\rho = -q(1+\psi_\rho).
\end{equation}
The space $L^2_\delta$ for $\delta\in\Rm$ is defined as the completion of $C^\infty_0(\Rm^n)$ with respect to the norm $\|\cdot\|_{L^2_\delta}$
defined as
\begin{equation}
  \label{eq:normdelta}
  \|u\|_{L^2_\delta}=\Big(\dint_{\Rm^n}
     \aver{x}^{2\delta}|u|^2 dx\Big)^{\frac12}, \qquad
   \aver{x}=(1+|x|^2)^{\frac12}.
\end{equation}
Let $-1<\delta<0$ and $q\in L^2_{\delta+1}$ and $\aver{x}q\in L^\infty$. One of the main results in \cite{Syl-Uhl-87} is that there exists  $\eta=\eta(\delta)$ such that the above problem admits a 
unique solution with $\psi_\rho\in L^2_\delta$ provided that 
\begin{displaymath}
  \|\aver{x}q\|_{L^\infty} +1 \leq \eta |\rho|.
\end{displaymath}
Moreover, $\|\psi_\rho\|_{L^2_\delta}\leq C|\rho|^{-1}\|q\|_{L^2_{\delta+1}}$ for some $C=C(\delta)$. In the analysis of many hybrid problems, we need smoother CGO solutions than what was recalled above. We introduce the spaces $H^s_\delta$ for $s\geq0$ as the completion of $C^\infty_0(\Rm^n)$ with respect to the norm $\|\cdot\|_{H^s_\delta}$ defined as
\begin{equation}
  \label{eq:normHsdelta}
  \|u\|_{H^s_\delta}=\Big(\dint_{\Rm^n}
     \aver{x}^{2\delta}|(I-\Delta)^{\frac s2}u|^2 dx\Big)^{\frac12}.
\end{equation}
Here $(I-\Delta)^{\frac s2}u$ is defined as the inverse Fourier transform of $\aver{\xi}^s\hat u(\xi)$, where $\hat u(\xi)$ is the Fourier transform of $u(x)$. Then we have the following 
\begin{proposition}[\cite{BU-IP-10}]
  \label{prop:regul}
  Let $-1<\delta<0$ and $k\in\Nm^*$. Let $q\in H^{\frac n2+k+\eps}_{1}$ and
  hence in $ H^{\frac n2+k+\eps}_{\delta+1}$ and $\rho$ be such that
  \begin{equation} \label{eq:constq}
    \|q\|_{H_1^{\frac n2+k+\eps}} +1 \leq \eta |\rho|.
  \end{equation}
  Then $\psi_\rho$ the unique solution to \eqref{eq:psi} belongs
  to $H^{\frac n2+k+\eps}_\delta$ and 
  \begin{equation}
    \label{eq:cts}
    |\rho| \|\psi_\rho\|_{H^{\frac n2+k+\eps}_\delta}
   \leq C  \|q\|_{H^{\frac n2+k+\eps}_{\delta+1}}, 
  \end{equation}
  for a constant $C$ that depends on $\delta$ and $\eta$.
\end{proposition}
We also want to obtain estimates for $\psi_\rho$ and $u_\rho$
restricted to the bounded domain $X$. We have the following result.
\begin{corollary}[\cite{BU-IP-10}]
  \label{cor:regX}
  Let us assume the regularity hypotheses of the previous proposition.
  Then we find that
  \begin{equation}
    \label{eq:reg}
    |\rho| \|\psi_\rho\|_{H^{\frac n2+k+\eps}(X)}
    + \|\psi_\rho\|_{H^{\frac n2+k+1+\eps}(X)} 
   \leq C  \|q\|_{H^{\frac n2+k+\eps}(X)}.
  \end{equation}
\end{corollary}
These results show that for $\rho$ sufficiently large, $\psi_\rho$ is small compared to $1$ in the class $C^k(\bar X)$ by Sobolev imbedding.

Let $Y=H^p(X)$ and $\mathcal M$ the ball in $Y$ of functions with norm bounded by a fixed $M>0$. Not only do we have that $\psi_\rho$ is small for $|\rho|$ large, but we have the following Lipschitz stability with respect to changes in the potential $q(x)$:
\begin{lemma}[\cite{BRUZ-IP-11}]
  \label{lem:lipsch} Let $\psi_\rho$ be the solution of
  \begin{equation}
    \label{eq:psisigma}
    \Delta \psi_\rho + 2\rho\cdot \nabla \psi_\rho = -q(1+\psi_\rho),
  \end{equation}
  and $\tilde\psi_\rho$ be the solution of the same equation with $q$
  replaced by $\tilde q$, where $\tilde q$ is defined as in
  \eqref{eq:q} with $\sigma$ replaced by $\tilde\sigma$. We assume
  that $q$ and $\tilde q$ are in $\mathcal M$. Then there is a
  constant $C$ such that for all $\rho$ with $|\rho|\geq|\rho_0|$, we
  have
  \begin{equation}
    \label{eq:lipsch}
    \|\psi_\rho-\tilde\psi_\rho\|_Y \leq \dfrac{C}{|\rho|} \|\sigma-\tilde\sigma\|_Y.
  \end{equation}
\end{lemma}
This is the property used in \cite{BRUZ-IP-11} to show that $\conductivitypot$ in the TAT problem \eqref{eq:invbrdy}-\eqref{eq:IScID} solves the equation
\begin{displaymath} 
  \conductivitypot(x) = e^{(\rho+\bar\rho)\cdot x} H(x) - {\mathcal H}_f[\conductivitypot](x) \mbox{ on } X,
\end{displaymath}
where 
\begin{displaymath} 
{\mathcal H}_f[\conductivitypot](x) = \conductivitypot \big( \psi_f+\overline{\psi_f}+\psi_f\overline{\psi_f}(x)\big),
 \end{displaymath}
is a contraction map for $f$ in an open set of illuminations; see \cite{BRUZ-IP-11}. The result in Theorem \ref{thm:reconstTAT} then follows by a Banach fixed point argument.

\subsubsection{CGO solutions and elliptic equations.} 
Consider the more general elliptic equation
\begin{equation}\label{eq:ellCGO} 
  -\nabla\cdot\diffusion\nabla u+\absorption u=0 \mbox{ in } X,\qquad u=f\mbox{ on } \dX.
\end{equation}
Upon defining $v=\sqrt\diffusion u$, we find that 
\begin{displaymath}
   (\Delta+q) v = 0 \mbox{ in } X,\qquad q=-\dfrac{\Delta\sqrt\diffusion}{\sqrt\diffusion} -\dfrac{\absorption}{\diffusion}.
 \end{displaymath}
 In other words, we find CGO solutions for \eqref{eq:ellCGO} defined on $\Rm^n$ and of the form
 \begin{equation}
  \label{eq:CGOell}
  u_\rho(x) = \dfrac{1}{\sqrt\diffusion}e^{\rho\cdot x} \big(1+\psi_\rho(x)\big),
\end{equation}
with $|\rho|\psi_\rho(x)$ bounded uniformly provided that $\diffusion$ and $\absorption$ are sufficiently smooth coefficients.

\subsubsection{Application to qualitative properties of elliptic solutions}

\paragraph{Lower bound for the modulus of complex valued solutions.} 
The above results show that for $|\rho|$ sufficiently large, then $|u_\rho|$ is uniformly bounded from below by a positive constant on compact domains. Note that $u_\rho$ is complex valued and that its real and imaginary parts oscillate very rapidly. Indeed,
\begin{displaymath} 
   e^{\rho\cdot x} = e^{\rho_r\cdot x} \big(\cos(\rho_i\cdot x) + i \sin(\rho_i\cdot x)\big),
\end{displaymath}  
which is rapidly increasing in the direction $\rho_r$ and rapidly oscillating in the direction $\rho_i$. Nonetheless, on a compact domain such as $X$, then $|u_\rho|$ is uniformly bounded from below by a positive constant. 

Let now $f_\rho={u_{\rho}}_{|\dX}$ the trace of the CGO solution on $\dX$. Then for $f$ close to $f_\rho$ and $u$ the solution to, say, \eqref{eq:HelCGO} or \eqref{eq:ellCGO}, we also obtain that $|u|$ is bounded from below by a positive constant. Such results were used in \cite{T-IP-10}.

\paragraph{Lower bound for vector fields.} 
For vector fields, we have the following result
\begin{theorem} [\cite{BU-IP-10}]
  \label{thm:field}
  Let $u_{\rho_j}$ for $j=1,2$ be CGO solutions with $q$ as above for both $\rho_j$ and $k\geq1$ and with $c_0^{-1}|\rho_1|\leq |\rho_2|\leq c_0|\rho_1|$ for some $c_0>0$. Then we have
  \begin{equation}
    \label{eq:drift}
    \hat\beta:= \dfrac{1}{2|\rho_1|} e^{-(\rho_1+\rho_2)\cdot x}
   \Big( u_{\rho_1}\nabla u_{\rho_2} - u_{\rho_2}\nabla u_{\rho_1}\Big)
    = \dfrac{\rho_1-\rho_2}{2|\rho_1|}
     + \hat h,
  \end{equation}
  where the vector field $\hat h$ satisfies the constraint
  \begin{equation}
    \label{eq:consthrho}
    \|\hat h\|_{C^k(\bar X)} \leq \dfrac{C_0}{|\rho_1|},
  \end{equation}
  for some constant $C_0$ independent of $\rho_{j}$, $j=1,2$.
\end{theorem}
With $\rho_2=\overline{\rho_1}$ so that $u_{\rho_2}=\overline{u_{\rho_1}}$, the imaginary part of \eqref{eq:drift} is a vector field that does not vanish on $X$ for $|\rho_1|$ sufficiently large. Moreover, let $u_{\rho_1}=v+iw$ and $u_{\rho_2}=v-iw$ for $v$ and $w$ real-valued functions. Then the imaginary and real parts of \eqref{eq:drift} are given by 
\begin{displaymath} 
   \Im \hat\beta = \dfrac{1}{|\rho_1|}e^{-2\Re\rho_1 \cdot x} (w\nabla v - v\nabla w)= \dfrac{\Im\rho_1}{|\rho_1|}
     + \Im \hat h ,\qquad \Re\hat\beta=0.
 \end{displaymath}
Let $u_1$ and $u_2$ be solutions of the elliptic problem \eqref{eq:HelCGO} on $X$ such that $u_1+iu_2$ on $\dX$ is close to the trace of $u_{\rho_1}$. The above result shows that
\begin{displaymath} 
   |u_1\nabla u_2-u_2\nabla u_1| \geq c_0>0 \quad \mbox{ in } X.
 \end{displaymath}
This yields \eqref{eq:lowerbd} and the result on unique and stable reconstructions in QPAT.

The above derivation may be generalized to the vector field $\beta_\alpha$ in \eqref{eq:betaalpha} with applications in elastography. Indeed let us start from \eqref{eq:decurho} with $\rho=\bk+i\bl$ such that  $\bk\cdot\bl=0$ and $k:=|\bk|=|\bl|$. Then using Corollary \ref{cor:regX}, we find that the following holds
\begin{equation} \label{eq:gradCGO}
   \begin{array}{rclrcl}
   \Re u_\rho &=& e^{\bk\cdot x} (\fc + \varphi_\rho^r) , \qquad &\Im u_\rho &=& e^{\bk\cdot x} (\fs + \varphi_\rho^i)\\
   \nabla \Re u_\rho  &=&k e^{\bk\cdot x} (\fc \hat \bk - \fs \hat \bl + \chi^r_\rho) ,\qquad & \nabla \Im u_\rho  &=&k e^{\bk\cdot x} (\fs \hat \bk + \fc \hat \bl + \chi^i_\rho) 
   \end{array}
\end{equation}
where $\fc= \cos (\bl\cdot x)$, $\fs = \sin(\bl\cdot x)$, $\hat \bk=\frac{\bk}{|\bk|}$, $\hat \bl=\frac{\bl}{|\bk|}$ and where $|\rho||\zeta|$ is bounded as indicated in Corollary \ref{cor:regX} for $\zeta\in\{\varphi_\rho^r,\varphi_\rho^i,\chi^r_\rho,\chi^i_\rho\}$. 

Let $u_1$ on $\partial X$ be close to $\Re u_\rho$.  Then we find by continuity that $|\nabla u_1 |$ is close to $|\nabla \Re u_\rho|$ so that for $k$ sufficiently large, we find that
\begin{equation}
  \label{eq:lowernablau1}
  |\nabla u_1 | \geq c_0 > 0  \quad \mbox{ in } X.
\end{equation}
This proves that $H(x)=\conductivity|\nabla u|$ is bounded from below by a positive constant provided that the boundary condition $f$ is in a  well-chosen open set of illuminations.

For the application to elastography, define now
\begin{displaymath} 
  \beta_\alpha = \Im u_\rho\nabla \Re u_\rho  - \alpha \Re u_\rho \nabla \Im u_\rho, \quad \alpha>0.
\end{displaymath}
Then for $|\rho|>\rho_\alpha$ sufficiently large so that $|\zeta|<\frac{(\min(1,\alpha))^2}{4(1+\alpha)}$ for $\zeta\in\{\varphi_\rho^r,\varphi_\rho^i,\chi^r_\rho,\chi^i_\rho\}$, we verify using \eqref{eq:gradCGO} that 
\begin{equation}
  \label{eq:lowerbetaalpha}
  |\beta_\alpha| \geq k e^{2\bk\cdot x} \dfrac12 \Big(\big(\fc \fs(1-\alpha) \big)^2 + \big(\fs^2+\alpha \fc^2\big)\Big) \geq k e^{2\bk\cdot x} \dfrac12 (\min(1,\alpha))^2.
\end{equation}
This provides a lower bound for $\beta_\alpha$ uniformly on compact sets. For an open set of  illuminations $(f_1,f_2)$ close to the traces of $(\Im u_\rho,\Re u_\rho)$ on $\partial X$, we find by continuity that the vector field $\beta_\alpha=u_1\nabla u_2-\alpha u_2\nabla u_1$ in \eqref{eq:betaalpha} also has a norm bounded from below uniformly on $X$.

\paragraph{Lower bound for determinants.}
The reconstruction in Theorem \ref{thm:stab3d} requires that the determinants in \eqref{eq:condji} be bounded from below. In specific situations, for instance when the conductivity is close to a given constant, such a determinant is indeed bounded from below by a positive constant for a large class of boundary conditions. However, it has been shown in \cite{BMN-ARMA-04} that the determinant of the gradients of three solutions could change signs on a domain with conductivities with large gradient. Unlike what happens in two dimensions of space, it is therefore not possible in general to show that the determinant of gradients of solutions has a given sign. However, using CGO solutions, we can be assured that on given bounded domains, the larger of two determinants is indeed uniformly positive for well-chosen boundary conditions. 

Let $u_\rho(x)$ be given by \eqref{eq:CGOell} solution of the elliptic problem \eqref{eq:ellCGO}. Upon treating the term $\psi_\rho$ and its derivative as in \eqref{eq:gradCGO} above and making them arbitrary small by choosing $\rho$ sufficiently large, we find that $\sqrt\conductivity u_\rho=e^{\rho\cdot x}+ {\rm l.o.t.}$ so that to leading order,
$$
\sqrt\conductivity\nabla u_\rho= e^{\bk\cdot x}(\bk+i\bl)\big(\cos(\bl\cdot x)+i\sin(\bl\cdot x)\big) +\rm{ l.o.t. }, \quad \rho=\bk+i\bl.
$$
Let $n=3$ and $(e_1,e_2,e_3)$ a constant orthonormal frame of $\Rm^3$.
It remains to take the real and imaginary parts of the above terms and choose $\hat \bk=e_2$ or $\hat \bk=e_3$ with $\hat \bl=e_1$ to obtain, up to normalization and negligible contributions  (for $k=|\bk|$ sufficiently large), that for 
\begin{displaymath} 
\begin{array}{rclrcl}
   \tilde S_1 &=& e_2 \cos k x_1 - e_1 \sin  k x_1 \qquad &
   \tilde S_2 &=& e_1 \cos k x_1 + e_2 \sin  k x_1 \\[1mm]
   \tilde S_3 &=& e_3 \cos k x_1 - e_1 \sin  k x_1 \qquad &
   \tilde S_4 &=& e_2 \cos|k|x_1 + e_3 \sin |k|x_1,
\end{array}
\end{displaymath}
we verify that $\det(\tilde S_1,\tilde S_2,\tilde S_3)=-\cos k x_1$ and that $\det(\tilde S_1,\tilde S_2,\tilde S_4)=-\sin k x_1$. Upon changing the sign of $S_3$ or $S_4$ if necessary to make both determinants non-negative, we find that the maximum of these two determinants is always bounded from below by a positive constant uniformly on $X$. This result is sufficient to prove Theorem \ref{thm:stab3d}; see \cite{BBMT-11}.
 
\paragraph{Hyperbolicity of a Lorentzian metric.} As a final application of CGO solutions, we mention the proof that a given constant vector field remains a time-like vector of a Lorentzian metric. This finds applications in the proof of Theorem \ref{thm:LinearStab} in \cite{B-UMEIT-11}.

Indeed, let $\hat \bk$ be a given direction in $\Sm^{n-1}$ and $\rho=i\bk+\bk^\perp$ and $u_\rho=e^{\rho\cdot x}$, once again neglecting $\psi_\rho$.  The real and imaginary parts of $\nabla u_\rho$ are such that 
\begin{equation} \label{eq:theta2} 
   e^{-\bk^\perp\!\cdot x}\Im \nabla e^{\rho\cdot x} = |\bk|\theta(x),\,\,e^{-\bk^\perp\!\cdot x}\Re \nabla e^{\rho\cdot x}=|\bk|\theta^\perp(x),
 \end{equation}
where $\theta(x) = \hat \bk\cos \bk\!\cdot\! x + \hat \bk^\perp \sin \bk\!\cdot \!x$ and $\theta^\perp(x)=-\hat \bk\sin \bk\!\cdot \!x+\hat \bk^\perp\cos \bk\!\cdot\! x$. As usual, $\hat \bk=\frac{\bk}{|\bk|}$.

Define the Lorentzian metrics 
\begin{displaymath} 
   \mh_\theta=2\theta\otimes\theta-I, \quad \mh_{\theta^\perp}=2\theta^\perp\otimes\theta^\perp-I.
 \end{displaymath}
Note that $\theta(x)$ and $\theta^\perp(x)$ oscillate in the plane $(\bk,\bk^\perp)$. The given vector $\hat \bk$ thus cannot be a time like vector for one of the Lorentzian metrics for all $x\in X$ (unless $X$ is a domain included in a thin slab). However, in the vicinity of any point $x_0$, we can construct a linear combination $\psi(x)=\cos\alpha\,\theta(x)+\sin\alpha\,\theta^\perp(x)$ for $\alpha\in[0,2\pi)$ such that 
\begin{displaymath} 
   \mbox{$\hat \bk$ is a time-like vector for }\mh_\psi=2\psi\otimes\psi-I, \quad \mbox{ i.e., } \quad \mh_\psi(\hat \bk,\hat \bk) = 2(\psi\cdot\hat \bk)^2-1 >0,
 \end{displaymath}
 uniformly for $x$ close to $x_0$; see \cite{B-UMEIT-11} for more details. When $\theta(x)$ is constructed as $\widehat{\nabla u}$ for $u$ solution to \eqref{eq:HelCGO} or \eqref{eq:ellCGO} for boundary conditions $f$ close to the trace of the corresponding CGO solution $u_{\rho}$, then the Lorentzian metric $\mh_\psi$ constructed above still verifies that $\hat \bk$ is a time-like vector with $ \mh_\psi(\hat \bk,\hat \bk)$ uniformly bounded from below by a positive constant locally.
%

\section{Conclusions and perspectives}
\label{sec:conclu}        

Research in hybrid inverse problems has been very active in recent years, primarily in the mathematical and medical imaging communities but also in geophysical imaging, see e.g. \cite{W-SIAP-05} and references on the electro-kinetic effect. This review focused on time-independent equations primarily with scalar-valued solutions. We did not consider the body of work done in the setting of time-dependent measurements, which involves different techniques than those presented here; see e.g. \cite{MZM-IP-10} and references. We considered scalar equations with the exception of the system of Maxwell's equations as it appears in Thermo-Acoustic Tomography. Very few results exist for systems of equations. The diffusion and conductivity equations considered in this review involve a scalar coefficient $\conductivity$. The reconstruction of more general tensors remains an open problem.


Compared to boundary value inverse problems, inverse problems with internal measurements enjoy better stability estimates precisely because local information is available. However, the derivation of such stability estimates often requires that specific, qualitative properties of solutions be satisfied, such as for instance the absence of critical points. This imposes constraints on the illuminations (boundary conditions) used to generate the internal data that forms one of the most difficult mathematical questions raised by the hybrid inverse problems. 

What are the ``optimal'' illuminations (boundary conditions)  for a given class of unknown parameters and how robust will the reconstructions be when such illuminations are modified are questions that are not fully answered. The theory of complex geometrical optics (CGO) solutions provides a useful tool to address these questions and construct suitable illuminations or prove their existence in several cases of interest. Numerical simulations will presumably be of great help to better understand whether such theoretical predictions are useful or reasonable in practice. Many numerical simulations performed in two dimensions of space confirm the good stability properties predicted by theory \cite{ABCTF-SIAP-08,BR-IP-11,BRUZ-IP-11,CFGK-SJIS-09,GS-SIAP-09,KK-AET-11,NTT-IP-09}. The two dimensional setting is special as we saw in section \ref{sec:2d}.  Very few simulations have been performed in the theoretically more challenging case of three (or more) dimensions of space. Simulations in \cite{KK-AET-11} show very promising three dimensional reconstructions in the setting of diffusion coefficients that are close to the constant case, which is also understood theoretically since $|\nabla u|$ then does not vanish for a large class of boundary conditions.

The main interest of hybrid inverse problems is that they combine high contrast with high resolution. This translates mathematically into good (Lipschitz or H\"older) stability estimates. Ideally, we would like to reconstruct highly oscillatory coefficients with a minimal influence of the noise in the measurements. Yet, all the results presented in this review paper and the cited references require that the coefficients satisfy some unwanted smoothness properties. To focus on one example for concreteness, the reconstructions in Photo-Acoustic Tomography involve the solution of the transport equation \eqref{eq:trchi}, which is well-posed provided that the vector field $\beta$ is sufficiently smooth. Using theories of renormalization, the regularity of such vector fields can be decreased to $W^{1,1}$ or to the BV category \cite{A-IM-04,BC-JMA-06,diperna-lions}. Yet, $u_1\nabla u_2-u_2\nabla u_1$ is a priori only in $L^2$ when $\conductivity$ is arbitrary as a bounded coefficient \cite{H-AIHP-03}. The construction of CGO solutions presented in section \ref{sec:CGO} also requires sufficient smoothness of the coefficients. How such reconstructions and stability estimates might degrade in the presence of non-smooth coefficients is quite open. Note that many similar problems are also open for boundary-value inverse problems \cite{U-IP-09}.

Finally, we have assumed in this review that the first step of the hybrid inverse problems had been done accurately. In practice, this may not quite always be so. PAT and TAT require that we solve an inverse source problem for a wave equation, which is a difficult problem in the presence of partial data and variable sound speed and is not entirely understood when realistic absorbing effects are accounted for \cite{KS-ATT-11,SU-11}. In UMEIT and UMOT, we have assumed in the derivation in section \ref{sec:physmod} that standing plane waves could be generated. This is practically difficult to achieve and different (equivalent) mechanisms have been proposed \cite{ABCTF-SIAP-08,KK-AET-11}. In transient elastography, we have assumed that the full (scalar) displacement could be reconstructed as a function of time and space. This is also sometimes an idealized approximation of what can be achieved in practice \cite{MZM-IP-10}. Finally, we have assumed knowledge of the current $\conductivity|\nabla u|$ in CDII, which is also difficult to acquire in practical settings as typically only the $z$ component of the magnetic field $B_z$ can be constructed; see the recent review \cite{SW-SR-11}. The modeling of errors generated during the first step of the procedure and the influence that such errors may have on the reconstructions during the second step of the hybrid inverse problem remain active areas of research.


%
\section*{Acknowledgment} This review and several collaborative efforts that led to papers referenced in the review were initiated during the participation of the author to the program on Inverse Problems and Applications at  the Mathematical Sciences Research Institute, Berkeley, California, in the Fall of 2010. I would like to thank the organizers and in particular Gunther Uhlmann for creating a very stimulating research environment at MSRI. Partial funding of this work by the National Science Foundation is also greatly acknowledged. 

{\footnotesize

} 

\end{document}